\documentclass[a4paper,11pt]{article}
\usepackage{amssymb}

\newcommand{\R}{\mathbb{R}}

\newcommand{\N}{\mathbb{N}}

\newcommand{\beq}{\begin{equation} }
\newcommand{\eqq}{\end{equation} }
\newcommand{\cuad}{{\sqcap\kern-.68em\sqcup}}
\newcommand{\abs}[1]{\mid #1 \mid}
\newcommand{\norm}[1]{\|#1\|}

\newtheorem{definition}{Definition}[section]
\newtheorem{teo}{Theorem}[section]

\newtheorem{proposition}{Proposition}[section]

\newtheorem{lemma}{Lemma}[section]
\newtheorem{corollary}{Corollary}[section]
\newtheorem{remark}{Remark}[section]
\newcommand{\bremark}{\begin{remark} \em}
\newcommand{\eremark}{\end{remark} }

\def\beeq{\begin{equation}}
\def\eeq{\end{equation}}
\newcommand{\begeqaet}{\begin{eqnarray*}}
\newcommand{\eneqaet}{\end{eqnarray*}}

\begin{document}

\begin{center}{\bf  \Large   Semilinear fractional elliptic equations \qquad\medskip

 involving measures }\medskip
%%%%%%%%%%%%%%%%%%%%%%%%%%%%%%%%%%%%%%%%%%%%%%%%%%%%%%%%%%%%%%%%%%%%%%
%%%%%%%%%%%%%%%%%%%%%%%%%%%%%%%%%%%%%%%%%%%%%%%%%%%%%%%%%%%%%%%%%%%%%%
\medskip
\bigskip

{\bf Huyuan Chen\footnote{chenhuyuan@yeah.net} \quad \quad  Laurent
V\'{e}ron\footnote{Laurent.Veron@lmpt.univ-tours.fr}}

Laboratoire de Math\'{e}matiques et Physique Th\'{e}orique
\\ CNRS UMR 7350
\\  Universit\'{e} Fran\c{c}ois Rabelais, Tours, France\\[1mm]
\bigskip
\begin{abstract}
We study the existence of  weak solutions to (E) $ (-\Delta)^\alpha
u+g(u)=\nu $ in a bounded regular domain $\Omega$ in $\R^N (N\ge2)$
which vanish in $\R^N\setminus\Omega$, where $(-\Delta)^\alpha$
denotes the fractional Laplacian with $\alpha\in(0,1)$, $\nu$ is a
Radon measure and $g$ is a nondecreasing function satisfying some
extra hypotheses. When $g$ satisfies a subcritical integrability
condition, we prove the existence and  uniqueness of a weak solution
for problem (E) for any measure. In the case where  $\nu$ is Dirac
measure, we characterize the asymptotic behavior of the solution.
When $g(r)=|r|^{k-1}r$ with $k$ supercritical, we show that a
condition of absolute continuity of the measure with respect to some
Bessel capacity is a necessary and sufficient condition in order (E)
to be solved.
\end{abstract}
\end{center}
\tableofcontents \vspace{1mm}
  \vspace{1mm}


\vspace{2mm}
%%%%%%%%%%%%%%%%%%%%%%%%%%%%%%%%%%%%%%%%%%%%%%%%%%%%%%%%%%%%%%%%%%%%%%%%%%%%%%%%%%%%%%%%%%%%%%%%%%%%%%%%%%%%%%%%%%%%%%%%%%
%%%%%%%%%%%%%%%%%%%%%%%%%%%%%%%%%%%%%%%%%%%%%%%%%%%%%%%%%%%%%%%%%%%%%%%%%%%%%%%%%%%%%%%%%%%%%%%%%%%%%%%%%%%%%%%%%%%%%%%%%%

\setcounter{equation}{0}
\section{Introduction}
Let $\Omega\subset \R^N$ be an open bounded  $C^2$ domain and
$g:\R\mapsto\R$ be a continuous function. We are concerned with the
existence of weak solutions  to the semilinear fractional elliptic
problem
\begin{equation}\label{eq1.1}
 \arraycolsep=1pt
\begin{array}{lll}
 (-\Delta)^\alpha  u+g(u)=\nu\quad & \rm{in}\quad\Omega,\\[2mm]
 \phantom{   (-\Delta)^\alpha  +g(u)}
u=0\quad & \rm{in}\quad \Omega^c,
\end{array}
\end{equation}
where   $ \alpha\in(0,1)$, $\nu$ is a Radon measure such that
$\int_\Omega\rho^\beta d|\nu|<\infty$ for some $\beta\in[0,\alpha]$
and $\rho(x)=dist(x,\Omega^c)$. The fractional Laplacian
$(-\Delta)^\alpha $  is defined by
$$(-\Delta)^\alpha  u(x)=\lim_{\epsilon\to0^+} (-\Delta)_\epsilon^\alpha u(x),$$
where for $\epsilon>0$,
\begin{equation}\label{1.2}
(-\Delta)_\epsilon^\alpha  u(x)=-\int_{\R^N}\frac{ u(z)-
u(x)}{|z-x|^{N+2\alpha}}\chi_\epsilon(|x-z|) dz
\end{equation}
and
$$\chi_\epsilon(t)=\left\{ \arraycolsep=1pt
\begin{array}{lll}
0,\quad & \rm{if}\quad t\in[0,\epsilon],\\[2mm]
1,\quad & \rm{if}\quad t>\epsilon.
\end{array}
\right.$$ When $\alpha=1$, the semilinear elliptic problem
\begin{equation}\label{eq003}
 \arraycolsep=1pt
\begin{array}{lll}
 -\Delta  u+g(u)=\nu \quad & \rm{in}\quad\Omega,\\[2mm]
 \phantom{ -\Delta  +g(u)}
u=0  \quad & \rm{on}\quad \partial\Omega,
\end{array}
\end{equation}
 has been extensively studied by numerous authors in the
last 30 years. A fundamental contribution is due to Brezis
\cite{B12}, Benilan and Brezis \cite{BB11}, where $\nu$ is a bounded
measure in $\Omega$ and the function $g:\R\to\R$ is nondecreasing,
positive on $(0,+\infty)$ and satisfies the subcritical assumption:
$$\int_1^{+\infty}(g(s)-g(-s))s^{-2\frac{N-1}{N-2}}ds<+\infty.$$
They proved the existence and uniqueness of the solution for problem
(\ref{eq003}). Baras and Pierre \cite{BP2} studied (\ref{eq003})
when $g(u)=|u|^{p-1}u$ for $p>1$ and $\nu$ is absolutely continuous
with respect to the Bessel capacity $C_{2,\frac{p}{p-1}}$, to obtain
a solution. In \cite{V} V\'{e}ron  extended Benilan and Brezis
results in replacing the Laplacian by a general uniformly elliptic
second order differential operator with Lipschitz continuous
coefficients; he obtained existence and uniqueness results for
solutions, when $\nu\in\mathfrak{M}(\Omega,\rho^\beta)$ with
$\beta\in[0,1]$ where $\mathfrak{M}(\Omega,\rho^\beta)$ denotes the
space of Radon measures in $\Omega$ satisfying
\begin{equation}\label{radon}
\int_{\Omega}\rho^\beta d|\nu|<+\infty,
\end{equation}
$\mathfrak{M}(\Omega,\rho^0)=\mathfrak{M}^b(\Omega)$ is the set of
bounded Radon measures and $g$ is nondecreasing and satisfies the
$\beta$-subcritical assumption:
$$\int_1^{+\infty}(g(s)-g(-s))s^{-2\frac{N+\beta-1}{N+\beta-2}}ds<+\infty.$$
The study of general semilinear elliptic equations with measure data
have been investigated, such as the equations involving measures
boundary data which was initiated by Gmira and V\'{e}ron \cite{GV}
who adapted the method introduced by Benilan and Brezis to obtain
the existence and uniqueness of solution. This subject has been
vastly expanded in recent years, see the papers of Marcus and
V\'{e}ron \cite{MV1,MV2,MV3,MV4}, Bidaut-V\'{e}ron and Vivier
\cite{BV}, Bidaut-V\'{e}ron, Hung and V\'{e}ron \cite{Hung}.

Recently, great attention has been devoted to non-linear equations
involving fractional Laplacian or more general integro-differential
operators and we mention the reference \cite{CT,CS0,CS1,CFQ,CLO2,L,RS,S1}.
In particular, the authors in \cite{KPU} used the duality approach to study  the
equations of
$$(-\Delta)^\alpha v=\mu\quad{\rm{in}}\quad \R^N,$$
where $\mu$ is a Radon measure with compact support. In \cite{CFQ} the authors obtained the existence of
large solutions to equation
\begin{equation}\label{eq1.2}
 (-\Delta)^\alpha  u+g(u)=f  \quad  \rm{in}\quad\Omega,
\end{equation}
where $\Omega$ is a bounded regular domain. In \cite{CV1} we considered the properties of possibly singular solutions of
(\ref{eq1.2}) in punctured domain . It is a well-known fact
\cite{LV1} that for $\alpha=1$ the weak singular solutions of
(\ref{eq1.2}) in punctured domain are classified according the type
of singularities they admits: either weak singularities with Dirac
mass, or strong singularities which are the upper limit of solutions
with weak singularities. One of our interests is to extend these
properties to any $\alpha\in(0,1)$ and furthermore to consider
general Radon measures.

In this paper we study the existence and uniqueness of solutions of
(\ref{eq1.1}) in a measure framework. Before stating our main
theorem we make precise the notion of weak solution  used in this
article.
\begin{definition}\label{weak definition}
We say that $u$ is a weak solution of (\ref{eq1.1}), if $u\in
L^1(\Omega)$, $g(u)\in L^1(\Omega,\rho^\alpha dx)$  and
\begin{equation}\label{weak sense}
\int_\Omega [u(-\Delta)^\alpha\xi+g(u)\xi]dx=\int_\Omega\xi
d\nu,\quad \forall\xi\in \mathbb{X}_{\alpha},
\end{equation}
where $\mathbb{X}_{\alpha}\subset C(\R^N)$ is the space of functions
$\xi$ satisfying:\smallskip

\noindent (i) $\rm{supp}(\xi)\subset\bar\Omega$,\smallskip

\noindent(ii) $(-\Delta)^\alpha\xi(x)$ exists for all $x\in \Omega$
and $|(-\Delta)^\alpha\xi(x)|\leq C$ for some $C>0$,\smallskip

\noindent(iii) there exist $\varphi\in L^1(\Omega,\rho^\alpha dx)$
and $\epsilon_0>0$ such that $|(-\Delta)_\epsilon^\alpha\xi|\le
\varphi$ a.e. in $\Omega$, for all
$\epsilon\in(0,\epsilon_0]$.\smallskip
\end{definition}

We notice that
 for $\alpha=1$, the test space $\mathbb{X}_{\alpha}$ is
used as $C_0^{1,L}(\Omega)$, which has similar properties like $(i)$
and $(ii)$. The  counter part for the Laplacian of assumption
$(iii)$ would be that the difference quotient
$\nabla_{x_j,h}[u](.):=h^{-1}[\partial_{x_j} u(.+h{\bf
e}_j)-\partial_{x_j}u(.)]$ is bounded by an $L^1$-function, which is
true since
$$\nabla_{x_j,h}[u](x)=h^{-1}\int_{0}^{h}\partial^2_{x_j,x_j} u(x+s{\bf e}_j)ds.
$$

We denote by  $G$ the Green kernel of $(-\Delta)^\alpha$ in
$\Omega\time\Omega $ and by $\mathbb{G}[.]$ the Green operator
defined by
\begin{equation}\label{optimal0}
\mathbb{G}[f](x)=\int_{\Omega}G(x,y)f(y) dy,\qquad\forall f\in
L^1(\Omega,\rho^\alpha dx).
\end{equation}
For $N\geq 2$, $0<\alpha<1$ and $\beta\in[0,\alpha]$, we define the
critical exponent
\begin{equation}\label{optimal1} k_{\alpha,\beta}=\left\{
\arraycolsep=1pt
\begin{array}{lll}
\frac{N}{N-2\alpha},\quad &\rm{if}\quad
\beta\in[0,\frac{N-2\alpha}N\alpha],\\[2mm]
\frac{N+\alpha}{N-2\alpha+\beta},\quad &\rm{if}\quad
\beta\in(\frac{N-2\alpha}N\alpha,\alpha].
\end{array}
\right.
\end{equation}
Our main result is the following:
%%%%%%%%%%%%%%%%%%%%%%%%%%%%%%%%%%%%%%%%%%%%%%%%%%%%%%%%%%%%%%%%%%%%%%%%%%%%%%%%%%%%%%%%%%%%%%%%%%%%%%%%%%%%%%%%%%%%%%%%%%%%%%%%%%%%%%%%%%%%%%%%%%%%%%%%%%%%%%%%%%%%%%%%%%%%%%%%%%%%%%%%%%%%
\begin{teo}\label{teo 1}
Assume $\Omega\subset \R^N$ ($N\ge2$) is an open bounded $C^2$
domain, $\alpha\in(0,1)$, $\beta\in[0,\alpha]$ and
$k_{\alpha,\beta}$ is defined by (\ref{optimal1}). Let $g:\R\to\R$
be a continuous, nondecreasing function, satisfying
\begin{equation}\label{1.4}
g(r)r\ge0,\quad \forall r\in\R\quad {\rm{and}}\quad \int_1^\infty
(g(s)-g(-s))s^{-1-k_{\alpha,\beta}}ds<\infty.
\end{equation}
Then for any $\nu\in\mathfrak{M}(\Omega,\rho^\beta)$ problem
(\ref{eq1.1}) admits a unique weak solution $u$. Furthermore, the
mapping: $\nu\mapsto u$ is increasing and
\begin{equation}\label{1.5}
-\mathbb{G}(\nu_-)\le u\le \mathbb{G}(\nu_+)\quad \rm{a.e.\ in}\
\Omega
\end{equation}
where $\nu_+$ and $\nu_-$ are respectively the positive and negative
part in the Jordan decomposition  of $\nu$.
\end{teo}

We note that for $\alpha=1$ and $\beta\in[0,1)$, we have
\begin{equation}\label{1.3}
k_{1,\beta}>\frac{N+\beta}{N-2+\beta},
\end{equation}
where $k_{1,\beta}$ is given in (\ref{optimal1}) and the number in
right hand side of (\ref{1.3}) is from Theorem 3.7 in \cite{V}.
Inspired by \cite{GV,V},  the existence of solution could be
extended in assuming that
 $g:\Omega\times\R\to\R$
is continuous and satisfies {\it the
$(N,\alpha,\beta)$-weak-singularity assumption}, that is, there
exists $r_0>0$ such that
$$g(x,r)r\ge0,\quad \forall (x,r)\in\Omega\times(\R\setminus(-r_0,r_0)),$$
and
$$|g(x,r)|\le \tilde g(|r|),\quad \forall (x,r)\in\Omega\times\R,$$
where $\tilde g:[0,\infty)\to[0,\infty)$ is continuous,
nondecreasing and satisfies
$$
 \int_1^\infty \tilde
g(s)s^{-1-k_{\alpha,\beta}}ds<\infty.
$$

 {\it We also give a stability result which shows that problem (\ref{eq1.1}) is weakly closed in the space of measures  $\mathfrak{M}(\Omega,\rho^\beta)$}. In the last section we characterize the behaviour of the solution $u$ of (\ref{eq1.1}) when $\nu=\delta_a$ for some $a\in \Omega$. We also study the case where $g(r)=|r|^{k-1}r$ when $k\geq k_{\alpha,\beta}$, which doesn't satisfy (\ref{1.4}). {\it We show that a necessary and sufficient condition in order a weak solution to problem
 \begin{equation}\label{eq1.k}
 \arraycolsep=1pt
\begin{array}{lll}
 (-\Delta)^\alpha  u+|u|^{k-1}u=\nu\quad & \rm{in}\quad\Omega,\\[2mm]
 \phantom{   (-\Delta)^\alpha  +|u|^{k-1}u}
u=0\quad & \rm{in}\quad \Omega^c,
\end{array}
\end{equation}
to exist where  $\nu$ is a positive bounded measure is that $\nu$
vanishes on compact subsets $K$ of $\Omega$ with zero
$C_{2\alpha,k'}$ Bessel-capacity.}

The paper is organized as follows. In Section 2 we give some
properties of Marcinkiewicz spaces and obtain the optimal index $k$
for which there holds
\begin{equation}\label{1.3'}
\|\mathbb{G}(\nu)\|_{M^{k}(\Omega,\rho^\gamma  dx)}\le
C\|\nu\|_{\mathfrak{M}(\Omega,\rho^\beta)}.
\end{equation}
We also gives some integration by parts formulas and prove a Kato's
type inequalities. In Section  3, we prove Theorem \ref{teo 1}. It
Section 4 we give applications the cases where the measure is a
Dirac mass and   where the nonlinearity is a power function.

%%%%%%%%%%%%%%%%%%%%%%%%%%%%%%%%%%%%%%%%%%
%%%%%%%%%%%%%%%%%%%%%%%%%%%%%%%%%%%%%%%%%%%
%%%%%%%%%%%%%%%%%%%%%%%%%%%%%%%%%%%%%%%%%%%%%%
%%%%%%%%%%%%%%%%%%%%%%%%%%%%%%%%%%%%%%%%%%%%%%%%%%%%%%%%%%%%%%%%%%%%%%%%%%%%%%%%%%%%%%%%%%%%%%%%%%%%%%%%%%%%%%%%%%%%%%%%%%%%%%%%%%%%%%%%%%%%%
%%%%%%%%%%%%%%%%%%%%%%%%%%%%%%%%%%%%%%%%%%%%%%
%%%%%%%%%%%%%%%%LINEAR%ESTIMATE%%%%%%%%%%%%%%%%%%%%%%%%%%%%%%%%%%%%%%%%%%%%%%%%%%%%%%%%%%%%%%%%%%%%%%%%%%%%%%%%%%%%%%%%%%%%%%%%%%%%%%%%%%%%%%%
%%%%%%%%%%%%%%%%%%%%%%%%%%%%%%%%%%%%%%%%%%%%%%
%%%%%%%%%%%%%%%%%%%%%%%%%%%%%%%%%%%%%%%%%%%%%%%%%%%%%%%%%%%%%%%%%%%%%%%%%%%%%%%%%%%%%%%%%%%%%%%%%%%%%%%%%%%%%%%%%%%%%%%%%%%%%%%%%%%%%%%%%%%%%
%%%%%%%%%%%%%%%%%%%%%%%%%%%%%%%%%%%%%%%%%%%%%%
%%%%%%%%%%%%%%%%%%%%%%%%%%%%%%%%%%%%%%%%%%%%%%%%%%%%%%%%%%%%%%%%%%%%%%%%%%%%%%%%%%%%%%%%%%%%%%%%%%%%%%%%%%%%%%%%%%%%%%%%%%%%%%%%%%%%%%%%%%%%%

\setcounter{equation}{0}
\section{Linear estimates}

\subsection{The Marcinkiewicz spaces}

We recall the definition and basic properties of the Marcinkiewicz
spaces.

\begin{definition}
Let $\Omega\subset \R^N$ be an open domain and $\mu$ be a positive
Borel measure in $\Omega$. For $\kappa>1$,
$\kappa'=\kappa/(\kappa-1)$ and $u\in L^1_{loc}(\Omega,d\mu)$, we
set
\begin{equation}\label{mod M}
\|u\|_{M^\kappa(\Omega,d\mu)}=\inf\{c\in[0,\infty]:\int_E|u|d\mu\le
c\left(\int_Ed\mu\right)^{\frac1{\kappa'}},\ \forall E\subset\Omega\
\rm{Borel\ set}\}
\end{equation}
and
\begin{equation}\label{spa M}
M^\kappa(\Omega,d\mu)=\{u\in
L_{loc}^1(\Omega,d\mu):\|u\|_{M^\kappa(\Omega,d\mu)}<\infty\}.
\end{equation}
\end{definition}

$M^\kappa(\Omega,d\mu)$ is called the Marcinkiewicz space of
exponent $\kappa$ or weak $L^\kappa$ space and
$\|.\|_{M^\kappa(\Omega,d\mu)}$ is a quasi-norm. The following
property holds.

\begin{proposition}\label{pr 1} \cite{BBC,CC}
Assume  $1\le q< \kappa<\infty$ and $u\in L^1_{loc}(\Omega,d\mu)$.
Then there exists  $C(q,\kappa)>0$ such that
$$\int_E |u|^q d\mu\le C(q,\kappa)\|u\|_{M^\kappa(\Omega,d\mu)}\left(\int_E d\mu\right)^{1-q/\kappa},$$
for any Borel set $E$ of $\Omega$.
\end{proposition}
 %%%%%%%%%%%%%%%%%%%%%%%%%%%%%%%%%%%%%%%%%%%%%%%%%%%%%%%%%%%%%%%%%%%%%%%%%%%%%%%%%%%%%%%%%%%%

For $\alpha\in(0,1)$ and $\beta,\gamma\in[0,\alpha]$ we set
\begin{equation}\label{parameters}
k_1(t)=\frac\gamma\alpha+\frac{N-(N-2\alpha)\frac\gamma\alpha}{N-2\alpha+t},\quad
 k_2(t)=\gamma+
\frac{N-(N-2\alpha)\frac\gamma\alpha}{N-2\alpha+t}t
\end{equation}
and
\begin{equation}\label{t0}
t_{\alpha,\beta,\gamma}=\min\{t\in[0,\alpha]:\
\frac{k_2(t)}{k_1(t)}\ge \beta\}.
\end{equation}
 %%%%%%%%%%%%%%%%%%%%%%%%%%%%%%%%%%%%%%%%%%%%%%%%%%%%%%%%%%%%%%%%%%%%%%%%%%%%%%%%%%%%%%%%%%%%
 \begin{remark} The quantity
$t_{\alpha,\beta,\gamma}$ is well defined, since
\begin{eqnarray*}
\frac{k_2(\alpha)}{k_1(\alpha)}=\frac{\gamma+\alpha\frac{N-(N-2\alpha)\frac\gamma\alpha}{N-\alpha}}{\frac{\gamma}{\alpha}+\frac{N-(N-2\alpha)\frac\gamma\alpha}{N-\alpha}}=\alpha\ge\beta.
\end{eqnarray*}
\end{remark}
 %%%%%%%%%%%%%%%%%%%%%%%%%%%%%%%%%%%%%%%%%%%%%%%%%%%%%%%%%%%%%%%%%%%%%%%%%%%%%%%%%%%%%%%%%%%%

\begin{remark} The function
$t\mapsto k_1(t)$ is decreasing in $[0,\alpha]$ with the following
bounds
$$k_1(0)=\frac{N}{N-2\alpha} \;\mbox{ and }\;k_1(\alpha)=\frac{N+\gamma}{N-\alpha}>1.$$
\end{remark}
 %%%%%%%%%%%%%%%%%%%%%%%%%%%%%%%%%%%%%%%%%%%%%%%%%%%%%%%%%%%%%%%%%%%%%%%%%%%%%%%%%%%%%%%%%%%%

\begin{remark}
The function $t\mapsto\frac{k_2(t)}{k_1(t)}$ is increasing in
$[0,\alpha]$, since
\begin{eqnarray*}
\left(\frac{k_2(t)}{k_1(t)}\right)'=\frac{[N-(N-2\alpha)\frac\gamma\alpha](N+\gamma)}{k_1^2(t)}>0.
\end{eqnarray*}
As a consequence  (\ref{t0}) is equivalent to
\begin{equation}\label{t2}
t_{\alpha,\beta,\gamma}=\max\{0,t_\beta\},
\end{equation}
where
\begin{equation}\label{t1}
t_{\beta}=\frac{\beta
N-(N-2\alpha)\gamma}{N-(N-3\alpha+\beta)\frac{\gamma}{\alpha}}.
\end{equation}
 is the solution of $\frac{k_2(t)}{k_1(t)}=\beta$.
\end{remark}
 %%%%%%%%%%%%%%%%%%%%%%%%%%%%%%%%%%%%%%%%%%%%%%%%%%%%%%%%%%%%%%%%%%%%%%%%%%%%%%%%%%%%%%%%%%%%

\begin{proposition}\label{general}
Let  $\Omega\subset \R^N\ (N\ge2)$ be an open  bounded $C^2$ domain
and $\nu\in\mathfrak{M}(\Omega,\rho^\beta)$ with
$\beta\in[0,\alpha]$. Then
\begin{equation}\label{annex 0}
\|\mathbb{G}[\nu]\|_{M^{k_{\alpha,\beta,\gamma}}(\Omega,\rho^\gamma
dx)}\le C\|\nu\|_{\mathfrak{M}(\Omega,\rho^\beta)},
\end{equation}
where $\gamma\in[0,\alpha]$, $\mathbb{G}[\nu](x)=\int_\Omega G(x,y)d
\nu(y)$ where $G$ is Green's kernel of $(-\Delta)^\alpha$ and
\begin{equation}\label{annex 00}
k_{\alpha,\beta,\gamma}= \left\{ \arraycolsep=1pt
\begin{array}{lll}
\frac{N+\gamma}{N-2\alpha+\beta},\quad &\rm{if}\ \gamma\le
\frac{N\beta}{N-2\alpha},\\[2mm]
\frac{N}{N-2\alpha},\quad &\rm{if\ not}.
\end{array}
\right.
\end{equation}
\end{proposition}
%%%%%%%%%%%%PROOF%%%%%%%%%%%%%%%%%%%%%%%%%%%%%%%%%%%%%%%%%%%%%%%%%%%%%%%%%%%%%%%%%%%%%%%%%%%%%%%%%%%%%%%%%%%%%%%%%%%%%%%%%%%%%%%%%%%%%%%%%%%%%
{\bf Proof.}  For $\lambda>0$ and $y\in\Omega$, we denote
$$A_\lambda(y)=\{x\in\Omega\setminus\{y\}: G(x,y)>\lambda\}\ \ {\rm
{and}}\quad m_\lambda(y)=\int_{A_\lambda(y)}\rho^\gamma(x) dx.$$
From \cite{CS}, there exists $C>0$ such that for any $(x,y)\in
\Omega\times\Omega$, $x\neq y$,
\begin{equation}\label{annex 01}
G(x,y)\le C
\min\left\{\frac1{|x-y|^{N-2\alpha}},\frac{\rho^\alpha(x)}{|x-y|^{N-\alpha}},\frac{\rho^\alpha(y)}{|x-y|^{N-\alpha}}\right\}
\end{equation}
and
\begin{equation}\label{annex 010}
G(x,y)\le C \frac{\rho^\alpha(y)}{\rho^\alpha(x)|x-y|^{N-2\alpha}}.
\end{equation}
Therefore, if $\gamma\in[0,\alpha]$ and $x\in A_\lambda(y)$, there
holds
\begin{equation}\label{annex 1}
\rho^\gamma(x)\le
\frac{C\rho^\gamma(y)}{\lambda^{\frac\gamma\alpha}|x-y|^{(N-2\alpha)\frac\gamma\alpha}}.
\end{equation}
Let $t\in[0,\alpha]$ be such that $\frac{k_2(t)}{k_1(t)}\ge \beta$,
where $k_1(t)$ and $k_2(t)$ are given in (\ref{parameters}), then
$$G(x,y)\le\left(\frac{C}{|x-y|^{N-2\alpha}}\right)^{1-\frac
t\alpha}\left(\frac{C\rho^\alpha(y)}{|x-y|^{N-\alpha}}\right)^{\frac
t\alpha}= \frac{C \rho^t(y)}{|x-y|^{N-2\alpha+t}}.$$
 We observe that
\begin{eqnarray*}
A_\lambda(y)\subset \left\{x\in\Omega\setminus\{y\}: \frac{C
\rho(y)^t}{|x-y|^{N-2\alpha+t}}>\lambda\right\}\subset D_\lambda(y)
\end{eqnarray*}
where
$D_\lambda(y):=\left\{x\in\Omega:|x-y|<(\frac{C\rho^t(y)}{\lambda})^{\frac1{N-2\alpha+t}}\right\}$;
together with (\ref{annex 1}), this implies
\begin{eqnarray*}
m_\lambda(y)\le\int_{D_\lambda(y)}\frac{C\rho^\gamma(y)}{\lambda^{\frac\gamma\alpha}|x-y|^{(N-2\alpha)\frac\gamma\alpha}}dx
\le C\rho(y)^{k_2(t)} \lambda^{-k_1(t)}.
\end{eqnarray*}
%%%%%%%%%%%%%%%%%%%%%%%%%%%%%%%%%%%%%%%%%%
%%%%%%%%%%%%%%%%%%%%%%%%%%%%%%%%%%%%%%%%%%%
%%%%%%%%%%%%%%%%%%%%%%%%%%%%%%%%%%%%%%%%%%
%%%%%%%%%%%%%%%%%%%%%%%%%%%%%%%%%%%%%%%%%%%
%%%%%%%%%%%%%%%%%%%%%%%%%%%%%%%%%%%%%%%%%%
%%%%%%%%%%%%%%%%%%%%%%%%%%%%%%%%%%%%%%%%%%%
%%%%%%%%%%%%%%%%%%%%%%%%%%%%%%%%%%%%%%%%%%%
For any Borel set $E$ of $\Omega$, we have
\begin{eqnarray*}
\int_E G(x,y)\rho^{\gamma}(x)dx\le
\int_{A_\lambda(y)}G(x,y)\rho^{\gamma}(x)dx+\lambda\int_E
\rho^{\gamma}(x)dx
\end{eqnarray*}
and
\begin{eqnarray*}
\int_{A_\lambda(y)}G(x,y)\rho^{\gamma}(x)dx&=&-\int_\lambda^\infty s
dm_s(y)
\\&=&\lambda m_\lambda(y)+ \int_\lambda^\infty m_s(y)ds
\\&\le& C\rho(y)^{k_2(t)}\lambda^{1-k_1(t)}.
\end{eqnarray*}
Thus,
\begin{eqnarray*}
\int_E G(x,y)\rho^{\gamma}(x)dx\le
C\rho(y)^{k_2(t)}\lambda^{1-k_1(t)}+\lambda \int_E
\rho^{\gamma}(x)dx.
\end{eqnarray*}
By choosing $\lambda=[\rho(y)^{-k_2(t)}\int_E
\rho^{\gamma}(x)dx]^{-\frac1{k_1(t)}}$, we have
\begin{eqnarray*}
\int_E G(x,y)\rho^{\gamma}(x)dx\le C\rho(y)^{\frac{k_2(t)}{k_1(t)}}
(\int_E \rho^{\gamma}(x)dx)^{\frac{k_1(t)-1}{k_1(t)}}.
\end{eqnarray*}
Therefore,
\begin{eqnarray*}
 \int_E
\mathbb{G}(|\nu|)(x)\rho^\gamma(x)dx&=&\int_\Omega\int_E
G(x,y)\rho^\gamma(x)dx d|\nu(y)|
\\&\le &C\int_\Omega \rho(y)^{\frac{k_2(t)}{k_1(t)}}d|\nu(y)|\left(\int_E
\rho^{\gamma}(x)dx\right)^{\frac{k_1(t)-1}{k_1(t)}}
\\&\le& C\|\nu\|_{\mathfrak{M}(\Omega,\rho^\beta)} \left(\int_E
\rho^{\gamma}(x)dx\right)^{\frac{k_1(t)-1}{k_1(t)}},
\end{eqnarray*}
since by our choice of $t$, $\frac{k_2(t)}{k_1(t)}\ge \beta$, which
guarantees that
  $$\int_\Omega
\rho(y)^{\frac{k_2(t)}{k_1(t)}}d|\nu(y)|\le
\max_{\Omega}\rho^{\frac{k_2(t)}{k_1(t)}-\beta}\int_\Omega
\rho(y)^{\beta}d|\nu(y)|.$$ As a consequence,
\begin{eqnarray*}
\|\mathbb{G}(\nu)\|_{M^{k_1(t)}(\Omega,\rho^\gamma dx)}\le
C\|\nu\|_{\mathfrak{M}(\Omega,\rho^\beta)}.
\end{eqnarray*}
Therefore,
$$k_{\alpha,\beta,\gamma}:=\max\{k_1(t):t\in [0,\alpha]\}=k_1(t_{\alpha,\beta,
\gamma}),
$$
where $t_{\alpha,\beta, \gamma}$ is defined by (\ref{t0}) and
$k_{\alpha,\beta,\gamma}$ is given by (\ref{annex 00}). We complete
the proof. \hfill$\Box$
%%%%%%%%%%%%%%%%%%%%%%%%%%%%%%%%%%%%%%%%%%
%%%%%%%%%%%%%%%%%%%%%%%%%%%%%%%%%%%%%%%%%%%
%%%%%%%%%%%%%%%%%%%%%%%%%%%%%%%%%%%%%%%%%%
%%%%%%%%%%%%%%%%%%%%%%%%%%%%%%%%%%%%%%%%%%%

\vspace{2mm}

We choose the parameter $\gamma$ in order to make
$k_{\alpha,\beta,\gamma}$ the largest possible, and denote
\begin{equation}\label{power1}
k_{\alpha,\beta}=\max_{\gamma\in[0,\alpha]} k_{\alpha,\beta,\gamma}.
\end{equation}
Since $\gamma\mapsto k_{\alpha,\beta,\gamma}$ is increasing, the
following statement holds.
\begin{proposition}\label{pr 3}
Let $N\ge2$ and $k_{\alpha,\beta}$ be defined by (\ref{power1}),
then
 \begin{equation}\label{power}
k_{\alpha,\beta}=\left\{ \arraycolsep=1pt
\begin{array}{lll}
\frac{N}{N-2\alpha},\quad &\rm{if}\quad
\beta\in[0,\frac{N-2\alpha}N\alpha],\\[2mm]
\frac{N+\alpha}{N-2\alpha+\beta},\quad &\rm{if}\quad
\beta\in(\frac{N-2\alpha}N\alpha,\alpha].
\end{array}
\right.
\end{equation}
\end{proposition}
%%%%%%%%%%%%%%%%%%%%%%%%%%%%%%%%%%%%%%%%%%
%%%%%%%%%%%%%%%%%%%%%%%%%%%%%%%%%%%%%%%%%%%
%%%%%%%%%%%%%%%%%%%%%%%%%%%%%%%%%%%%%%%%%%%%%%
%%%%%%%%%%%%%%%%%%%%%%%%%%%%%%%%%%%%%%%%%%%%%%%%%%%%%%%%%%%%%%%%%%%%%%%%%%%%%%%%%%%%%%%%%%%%%%%%%%%%%%%%%%%%%%%%%%%%%%%%%%%%%%%%%%%%%%%%%%%%%
%%%%%%%%%%%%%%%%%%%%%%%%%%%%%%%%%%%%%%%%%%%%%%
%%%%%%%%%%%%%%%%%%%%%%%%%%%%%%%%%%%%%%%%%%%%%%%%%%%%%%%%%%%%%%%%%%%%%%%%%%%%%%%%%%%%%%%%%%%%%%%%%%%%%%%%%%%%%%%%%%%%%%%%%%%%%%%%%%%%%%%%%%%%%
%%%%%%%%%%%%%%%%%%%%%%%%%%%%%%%%%%%%%%%%%%%%%%
%%%%%%%%%%%%%%%%%%%%%%%%%%%%%%%%%%%%%%%%%%%%%%%%%%%%%%%%%%%%%%%%%%%%%%%%%%%%%%%%%%%%%%%%%%%%%%%%%%%%%%%%%%%%%%%%%%%%%%%%%%%%%%%%%%%%%%%%%%%%%
%%%%%%%%%%%%%%%%%%%%%%%%%%%%%%%%%%%%%%%%%%%%%%
%%%%%%%%%%%%%%%%%%%%%%%%%%%%%%%%%%%%%%%%%%%%%%%%%%%%%%%%%%%%%%%%%%%%%%%%%%%%%%%%%%%%%%%%%%%%%%%%%%%%%%%%%%%%%%%%%%%%%%%%%%%%%%%%%%%%%%%%%%%%%
%%%%%%%%%%%%%%%%%%%%%%%%%%%%%%%%%%%%%%%%%%%
%%%%%%%%%%%%%%%%NON-HOMOGENEOUS PROBLEMS%%%%%%%%%%%%%%
%%%%%%%%%%%%%%%%%%%%%%%%%%%%%%%%%%%%%%%%%%%
%%%%%%%%%%%%%%%%%%%%%%%%%%%%%%%%%%%%%%%%%%
%%%%%%%%%%%%%%%%%%%%%%%%%%%%%%%%%%%%%%%%%%%
%%%%%%%%%%%%%%%%%%%%%%%%%%%%%%%%%%%%%%%%%%%

\subsection{ Non-homogeneous problem}

In this subsection, we study some properties of the solution of the
linear non-homogeneous, which will play a key role in the sequel. We
assume that $\Omega\subset\R^N$, $N\geq 2$ is an open bounded domain
with a $C^2$ boundary.

\begin{lemma}\label{eta}
$(i)$\ There exists $C>0$ such that for any $\xi\in
\mathbb{X}_{\alpha}$ there holds
%$$|\xi(x)|\le C\norm{(-\Delta)^\alpha\xi}_{L^\infty(\Omega)}\rho^\alpha(x),\quad x\in\Omega,$$
\begin{equation}\label{est 0}
\norm{\xi}_{C^\alpha(\bar\Omega)}\le C
\norm{(-\Delta)^\alpha\xi}_{L^\infty(\Omega)}
\end{equation}
and
\begin{equation}\label{est 1}
\norm{\rho^{-\alpha}\xi}_{C^\theta(\bar\Omega)}\le C
\norm{(-\Delta)^\alpha\xi}_{L^\infty(\Omega)}.
\end{equation}
where $0<\theta<\min\{\alpha, 1-\alpha\}$. In particular, for
$x\in\Omega$
\begin{equation}\label{est 2}
|\xi(x)|\le
C\norm{(-\Delta)^\alpha\xi}_{L^\infty(\Omega)}\rho^\alpha(x).
\end{equation}

\noindent$(ii)$ Let $u$ be the solution of
\begin{equation}\label{eta 1}
 \arraycolsep=1pt
\begin{array}{lll}
 (-\Delta)^\alpha  u=f\quad & \rm{in}\quad\Omega,\\[2mm]
 \phantom{ (-\Delta)^\alpha}
u=0\quad & \rm{in}\quad \Omega^c,
\end{array}
\end{equation}
where $f\in C^\gamma(\bar \Omega)$ for $\gamma>0$.
 Then $u\in \mathbb{X}_{\alpha}$.
\end{lemma}
%%%%%%%%%%%%PROOF%%%%%%%%%%%%%%%%%%%%%%%%%%%%%%%%%%%%%%%%%%%%%%%%%%%%%%%%%%%%%%%%%%%%%%%%%%%%%%%%%%%%%%%%%%%%%%%%%%%%%%%%%%%%%%%%%%%%%%%%%%%%%
{\bf Proof.} $\textbf{(i)}$. Estimates $(\ref{est 0})$ and
$(\ref{est 1})$ are consequences of \cite[Prop 1.1]{RS} and \cite[Th
1.2]{RS} respectively. Furthermore, if $\eta_1$ is the solution of
(\ref{eta 1}) with $f\equiv1$ in $\Omega$,  then $\eta_1>0$ in
$\Omega$ and by follows \cite[Th 1.2]{RS}, there exists $C>0$ such
that
\begin{equation}\label{est1-1}
C^{-1}\leq
\frac{\eta_1}{\rho^\alpha}\le C\quad\rm{in}\quad \Omega.
\end{equation}
In this expression the right-side follows \cite[Th 1.2]{RS}
 and  the left-hand side inequality follows from the maximum principle and \cite[Th 1.2]{CS}.
Since
$$-\norm{(-\Delta)^\alpha\xi}_{L^\infty(\Omega)}\le(-\Delta)^\alpha\xi\le
\norm{(-\Delta)^\alpha\xi}_{L^\infty(\Omega)}\quad \rm{in}\
\Omega,$$ it follows by the comparison principle,
$$
-\norm{(-\Delta)^\alpha\xi}_{L^\infty(\Omega)}\eta_1(x)\le\xi(x)\le
\norm{(-\Delta)^\alpha\xi}_{L^\infty(\Omega)}\eta_1(x).
$$
which, together with (\ref{est1-1}), implies (\ref{est 2}).

\noindent$\textbf{(ii)}$ For $r>0$, we denote
$$\Omega_{r}=\{z\in\Omega:\ dist(z,\partial\Omega)>r\}.$$ Since
$f\in C^\gamma(\bar \Omega)$, then by  Corollary 1.6 part $(i)$ and
Proposition 1.1  in \cite{RS}, for
$\theta\in[0,\min\{\alpha,1-\alpha,\gamma\})$, there exists $C>0$
such that for any $r>0$, we have
$$\|u\|_{C^{2\alpha+\theta}(\Omega_r)}\le Cr^{-\alpha-\theta}$$
and
$$\|u\|_{C^{\alpha}(\R^N)}\le C.$$
Then for $x\in\Omega$, letting $r=\rho(x)/2$,
\begin{equation}\label{wy1}
|\delta(u,x,y)|\le Cr^{-\alpha-\theta}|y|^{2\alpha+\theta},\quad
\forall y\in B_r(0)
\end{equation}
 and
$$|\delta(u,x,y)|\le C|y|^{\alpha},\quad \forall y\in \R^N,$$
where $\delta(u,x,y)=u(x+y)+u(x-y)-2u(x)$. Thus,
\begin{eqnarray*}
 |(-\Delta)_\epsilon^\alpha
 u(x)|&\le&\frac12\int_{\R^N}\frac{|\delta(u,x,y)|}{|y|^{N+2\alpha}}\chi_\epsilon
 (|y|)dy
 \\&\le&\frac12\int_{B_{r}(0)}\frac{|\delta(u,x,y)|}{|y|^{N+2\alpha}}dy+\frac12\int_{B^c_{r}(0)}\frac{|\delta(u,x,y)|}{|y|^{N+2\alpha}}dy
 \\&\le&\frac{Cr^{-\alpha-\theta}}2\int_{B_r(0)}\frac1{|y|^{N-\theta}}dy+\frac C2\int_{B^c_r(0)}\frac1{|y|^{N+\alpha}}dy
 \\&\le& C\rho(x)^{-\alpha},\quad x\in\Omega,
\end{eqnarray*}
for some  $C>0$ independent of $\epsilon$. Moreover,
$\rho^{-\alpha}$ is
 in $%L^1(\Omega)\subset
 L^1(\Omega,\rho^\alpha dx)$.
Finally, we prove $(-\Delta)_\epsilon^\alpha
 u\to(-\Delta)^\alpha
 u$ as $\epsilon\to0^+$ pointwise.
For  $x\in\Omega$, choosing $\epsilon\in(0,\rho(x)/2)$, then by
(\ref{wy1}),
\begin{eqnarray*}
 |(-\Delta)^\alpha
 u(x)-(-\Delta)_\epsilon^\alpha
 u(x)|&\le&\frac12\int_{B_\epsilon(0)}\frac{|\delta(u,x,y)|}{|y|^{N+2\alpha}}dy
 %\\&\le&\frac{Cd(x)^{-\alpha-\theta}}2\int_{B_\epsilon(0)}\frac1{|y|^{N+2\alpha-2}}dy
 \\&\le& C\rho(x)^{-\alpha-\theta}\epsilon^\theta
 \\&\to&0,\quad \epsilon\to0^+.
\end{eqnarray*}
 The proof is
complete. \hfill$\Box$

%\begin{definition}
%We say $u$ is a weak solution of
%\begin{equation}\label{homo}
% \arraycolsep=1pt
%\begin{array}{lll}
 %(-\Delta)^\alpha u=\nu,\quad & \rm{in}\quad\Omega,\\[2mm]
%u=0,\quad & \rm{in}\quad \R^N\setminus\Omega,
%\end{array}
%
%\end{equation}
%where $\nu\in\mathfrak{M}(\Omega,\rho^\alpha dx)$, if   $u\in
%L^1(\Omega)$ and
%\begin{equation}\label{weak dh}
%\int_\Omega u (-\Delta)^\alpha \xi dx=\int_\Omega \xi d\nu,\quad
%\forall\xi\in\mathbb{X}_{\alpha}.
%\end{equation}
%\end{definition}
%It is well defined for (\ref{weak dh}) by Lemma \ref{eta}.
\medskip

The following Proposition is the Kato's type estimate for proving
the uniqueness of the solution of (\ref{eq1.1}).

\begin{proposition}\label{pr 2.1}
If $\nu\in L^1(\Omega,\rho^\alpha dx)$, there exists a  unique weak
solution $u$ of the problem
\begin{equation}\label{homo}
 \arraycolsep=1pt
\begin{array}{lll}
 (-\Delta)^\alpha &u=\nu\quad & \rm{in}\quad\Omega,\\[2mm]
&u=0\quad & \rm{in}\quad \Omega^c.
\end{array}
\end{equation}
For any $\xi\in\mathbb{X}_{\alpha}$, $\xi\ge0$, we have
 \begin{equation}\label{sign}
\int_\Omega |u|(-\Delta)^\alpha \xi dx\le \int_\Omega  \xi
{\rm{sign}}(u)\nu dx
\end{equation}
and
 \begin{equation}\label{sign+}
\int_\Omega u_+(-\Delta)^\alpha \xi dx\le \int_\Omega  \xi
{\rm{sign}}_+(u)\nu dx,
\end{equation}
%and
%\begin{equation}\label{sign0}
%\int_\Omega u_+(-\Delta)^\alpha \xi dx\le \int_\Omega  \xi
%\rm{sign}_+(u) \nu dx,
%\end{equation}
%where $u_+=\max\{0,u\}$, $\rm{sign}_+(u)=1$ if $u>0$ and
%$\rm{sign}_+(u)=0$ if $u\le0$.
\end{proposition}

We note here that for $\alpha=1$, the proof of Proposition \ref{pr
2.1} could be seen in \cite[Th 2.4]{V}.  For $\alpha\in(0,1)$, we
first prove some integration by parts formula.
%%%%%%%%%%%%%%%%%%%%%%%%%%%%%%%%%%%%%%%%%%%%%%%%%%%%%%%%%%%%%%%%%%%%%%%%%%%%%%%%%%%%%%%%%%%%%%%%%%%%%%%%%%%%%%%%%%%%%%%%%%%%%%%%%%%%%%%%%%%%

\begin{lemma}\label{lm 2.1}
Assume $u,\xi\in \mathbb{X}_{\alpha}$, then
\begin{equation}\label{W n}
\int_\Omega u (-\Delta)^\alpha\xi dx=\int_\Omega\xi (-\Delta)^\alpha
u dx.
\end{equation}
\end{lemma}
%%%%%%%%%%%%PROOF%%%%%%%%%%%%%%%%%%%%%%%%%%%%%%%%%%%%%%%%%%%%%%%%%%%%%%%%%%%%%%%%%%%%%%%%%%%%%%%%%%%%%%%%%%%%%%%%%%%%%%%%%%%%%%%%%%%%%%%%%%%%%
{\bf Proof.} Denote
 \begin{equation}\label{4.2}
 (-\Delta)_{\Omega,\epsilon}^\alpha
u(x)=-\int_{\Omega}\frac{u(z)-u(x)}{|z-x|^{N+2\alpha}}\chi_\epsilon(|x-z|)dz.
\end{equation}
By the definition of $(-\Delta)_\epsilon^\alpha$, we have
\begin{eqnarray*}
 (-\Delta)_\epsilon^\alpha
 u(x)&=&-\int_{\Omega}\frac{u(z)-u(x)}{|z-x|^{N+2\alpha}}\chi_\epsilon(|x-z|)dz+u(x)\int_{\Omega^c}\frac{\chi_\epsilon(|x-z|)}{|z-x|^{N+2\alpha}}dz
 \\&=& (-\Delta)_{\Omega,\epsilon}^\alpha
 u(x)+u(x)\int_{\Omega^c}\frac{\chi_\epsilon(|x-z|)}{|z-x|^{N+2\alpha}}dz.
\end{eqnarray*}
We claim that
\begin{equation}\label{2.1}
\int_{\Omega}\xi(x)(-\Delta)_{\Omega,\epsilon}^\alpha
u(x)dx=\int_{\Omega}u(x)(-\Delta)_{\Omega,\epsilon}^\alpha \xi(x)
dx,\quad {\rm{for}}\  u,\xi \in\mathbb{X}_{\alpha}.
\end{equation}
By using the fact of
\begin{eqnarray*}
\int_{\Omega}\int_{\Omega}\frac{[u(z)-u(x)]\xi(x)}{|z-x|^{N+2\alpha}}\chi_\epsilon(|x-z|)dzdx
=\int_{\Omega}\int_{\Omega}\frac{[u(x)-u(z)]\xi(z)}{|z-x|^{N+2\alpha}}\chi_\epsilon(|x-z|)dzdx,
\end{eqnarray*}
 we have
\begin{eqnarray*}
&&\int_{\Omega}\xi(x)(-\Delta)_{\Omega,\epsilon}^\alpha
u(x)dx\\&&=-\frac{1}2\int_{\Omega}\int_{\Omega}[\frac{(u(z)-u(x))\xi(x)}{|z-x|^{N+2\alpha}}+\frac{(u(x)-u(z))\xi(z)}{|z-x|^{N+2\alpha}}]\chi_\epsilon(|x-z|)dzdx
\\&&=\frac{1}2\int_{\Omega}\int_{\Omega}\frac{[u(z)-u(x)][\xi(z)-\xi(x)]}{|z-x|^{N+2\alpha}}\chi_\epsilon(|x-z|)dzdx.
\end{eqnarray*}
Similarly, by the fact that $u\in\mathbb{X}_{\alpha}$,
\begin{eqnarray*}
\int_{\Omega}u(x)(-\Delta)_{\Omega,\epsilon}^\alpha \xi(x)
dx=\frac12\int_{\Omega}\int_{\Omega}\frac{[u(z)-u(x)][\xi(z)-\xi(x)]}{|z-x|^{N+2\alpha}}\chi_\epsilon(|x-z|)dzdx.
\end{eqnarray*}
Then (\ref{2.1}) holds. In order to prove (\ref{W n}), we first
notice that by (\ref{2.1}),
\begin{eqnarray}
&&\int_{\Omega}\xi(x)(-\Delta)_\epsilon^\alpha u(x)dx \nonumber
\\&&=\int_{\Omega}\xi(x)(-\Delta)_{\Omega,\epsilon}^\alpha u(x)dx+\int_\Omega
u(x)\xi(x)\int_{\Omega^c}\frac{\chi_\epsilon(|x-z|)}{|z-x|^{N+2\alpha}}dzdx\nonumber
\\&&=\int_{\Omega} u(x)(-\Delta)_{\Omega,\epsilon}^\alpha\xi(x)dx+\int_\Omega
u(x)\xi(x)\int_{\Omega^c}\frac{\chi_\epsilon(|x-z|)}{|z-x|^{N+2\alpha}}dzdx\nonumber
\\&&=\int_{\Omega}u(x)(-\Delta)_\epsilon^\alpha
\xi(x)dx.\label{2.4}
\end{eqnarray}
Since  $u$ and $\xi$ belongs to $\mathbb{X}_{\alpha}$,
$(-\Delta)_\epsilon^\alpha \xi\to (-\Delta)^\alpha \xi$ and
$(-\Delta)_\epsilon^\alpha u\to (-\Delta)^\alpha u$ and
$|u(-\Delta)_\epsilon^\alpha \xi|+|\xi(-\Delta)_\epsilon^\alpha
u|\leq C\varphi$ for some $C>0$ and $\varphi\in
L^1(\Omega,\rho^\alpha dx)$. It follows by the Dominated Convergence
Theorem
$$\lim_{\epsilon\to0^+}\int_{\Omega}\xi(x)(-\Delta)_\epsilon^\alpha
u(x)dx=\int_{\Omega}\xi(x)(-\Delta)^\alpha u(x)dx$$ and $$
\lim_{\epsilon\to0^+}\int_{\Omega}(-\Delta)_\epsilon^\alpha\xi(x)
u(x)dx=\int_{\Omega}(-\Delta)^\alpha \xi(x)u(x)dx.$$ Letting
$\epsilon\to0^+$ of (\ref{2.4}) we conclude that (\ref{W n}) holds.
\hfill$\Box$

\medskip

%%%%%%%%%%%%%%%%%%%%%%%%%%%%%%%%%%%%%%%%%%%%%%%%%%%%%%%%%%%%%%%%%%%%%%%%%%%%%%%%%%%%%%%%%%%%%%%%%%%%%%%%%%%%%%%%%%%%%%%%%%%%%%%%%%%%%%%%%%%%%%%%%%%%%%%PROOF%%%%%%%%%%%%%%%%%%%%%%%%%%%%%%%%%%%%%%%%%%%%%%%%%%%%%%%%%%%%%%%%%%%%%%%%%%%%%%%%%%%%%%%%%%%%%%%%%%%%%%%%%%%%%%%%%%%%%%%%%%%%%%

For $1\leq p<\infty$ and $0<s<1$,  $W^{s,p}(\Omega)$ is the set of
$\xi\in L^p(\Omega)$ such that
 \begin{equation}\label {X2}
\int_{\Omega}\int_{\Omega}\frac{|\xi(x)-\xi(y)|^p}{|x-y|^{N+s
p}}dydx<\infty.
\end{equation}
This space is endowed with the norm
\begin{equation}\label {X3}
\norm\xi_{W^{s,p}(\Omega)}=\left(\int_{\Omega}|\xi(x)|^p
dx+\int_{\Omega}\int_{\Omega}\frac{|\xi(x)-\xi(y)|^p}{|x-y|^{N+s
p}}dydx\right)^{\frac{1}{p}}.
\end{equation}
Furthermore, if $\Omega$ is bounded, the following Poincar\'e
inequality holds \cite[p 134]{T}.
\begin{equation}\label {X4}
\left(\int_{\Omega}|\xi(x)|^p dx\right)^{\frac{1}{p}}\leq
C\left(\int_{\Omega}\int_{\Omega}\frac{|\xi(x)-\xi(y)|^p}
{|x-y|^{N+s p}}dydx\right)^{\frac{1}{p}},\quad\forall\xi\in
C_c^{\infty}(\Omega).
\end{equation}
%Let $W^{s,p}_0(\Omega)$ denote the closure of $C^\infty_c(\Omega)$
%in the norm $\|\cdot\|_{W^{s,p}(\Omega)}$.

\begin{lemma}\label{lm 2.2}
  Let $u\in\mathbb{X}_{\alpha}$ and $\gamma$ be $C^2$ in the interval $u(\bar\Omega)$ and
satisfy $\gamma (0)=0$ , then $u\in W^{\alpha,2}(\Omega)$,
$\gamma\circ u\in\mathbb{X}_{\alpha}$ and for  all $x\in\Omega$,
there exists $z_x\in\bar\Omega$ such that
\begin{equation}\label {K1}
(-\Delta)^{\alpha}(\gamma\circ u)(x)=(\gamma'\circ
u)(x)(-\Delta)^{\alpha}u(x)- \frac{\gamma''\circ
u(z_x)}{2}\int_{\Omega}\frac{(u(y)-u(x))^2}{|y-x|^{N+2\alpha}}dy.
\end{equation}
\end{lemma}
%%%%%%%%%%%%PROOF%%%%%%%%%%%%%%%%%%%%%%%%%%%%%%%%%%%%%%%%%%%%%%%%%%%%%%%%%%%%%%%%%%%%%%%%%%%%%%%%%%%%%%%%%%%%%%%%%%%%%%%%%%%%%%%%%%%%%%%%%%%%%
{\bf Proof.}  Since $u\in C(\bar\Omega)$ vanishes in $\Omega^c$,
$\gamma\circ u$ shares the same properties. By $(\ref{est 0})$, for
any $x$ and $y$ in $\Omega$
$$(u(x)-u(y))^2\leq
C|x-y|^{2\alpha}\norm{(-\Delta)^{\alpha}u}^2_{L^\infty(\Omega)}.
$$
Then $u\in W^{\alpha,2}(\Omega)$. Similarly $\gamma\circ u\in
W^{\alpha,2}(\Omega)$. Furthermore
$$\begin{array} {ll}\displaystyle
(\gamma\circ u)(y)-(\gamma\circ u)(x)%\\[0mm]\phantom{(\gamma\circ u)(x+y))}\displaystyle
=(\gamma'\circ
u)(x)\left(u(y)-u(x)\right)+\int_{u(x)}^{u(y)}\!\!\!(u(y)-t)\gamma''(t)
dt.
\end{array}$$
By the mean value theorem, there exists some $\tau\in [0,1]$ such
that
$$\int_{u(x)}^{u(y)}\!\!\!(u(y)-t)\gamma''(t) dt=\frac{\gamma''(\tau
u(y)+(1-\tau)u(x))}{2}(u(y)-u(x))^2.
$$
Since $\gamma'' $ is continuous and $u$ is continuous in $\bar
\Omega$,
$$ \left|\int_{u(x)}^{u(y)}\!\!\!(u(y)-t)\gamma''(t) dt\right|\leq \frac{\norm{\gamma''\circ u}_{L^\infty(\bar\Omega)}}{2}(u(y)-u(x))^2
$$
and by $(\ref{est 0})$,
$$\begin{array} {ll}\displaystyle
\left|\int_{\abs {y-x}>\epsilon}
\int_{u(x)}^{u(y)}\!\!\!(u(y)-t)\gamma''(t) dt\frac{dy}{|y-x|^{N+2\alpha}}\right|\\[4mm]\displaystyle\phantom{--------}
\leq \frac{\norm{\gamma''\circ
u}_{L^\infty}}{2}\int_{\Omega}(u(y)-u(x))^2\frac{dy}{|y-x|^{N+2\alpha}}.
\end{array}$$
Notice also that $\tau u(y)+(1-\tau)u(x)\in u(\bar\Omega) :=I $,
therefore
$$\min_{t\in I} \gamma''(t)\leq\gamma''(\tau u(y)+(1-\tau)u(x))\leq \max_{t\in I}
\gamma''(t),
$$
thus
$$\begin{array} {ll}\displaystyle
\frac{\min_{t\in I} \gamma''(t)}{2}\int_{\Omega}
\frac{(u(y)-u(x))^2}{|y-x|^{N+2\alpha}}dy\\[4mm]\displaystyle\phantom{--------}
\leq \int_{\Omega}
\int_{u(x)}^{u(y)}\!\!\!(u(y)-t)\gamma''(t) dt \frac{dy}{|y-x|^{N+2\alpha}}\\[4mm]\displaystyle\phantom{--------------}\displaystyle\leq \frac{\max_{t\in I} \gamma''(t)}{2}\int_{\Omega}
\frac{(u(y)-u(x))^2}{|y-x|^{N+2\alpha}}dy.
\end{array}$$
Since $\gamma''$ is continuous, there exists $t_0\in I$ such that
$$\int_{\Omega}
\int_{u(x)}^{u(y)}\!\!\!(u(y)-t)\gamma''(t) dt
\frac{dy}{|y-x|^{N+2\alpha}}=\frac{\gamma''(t_0)}2\int_{\Omega}
\frac{(u(y)-u(x))^2}{|y-x|^{N+2\alpha}}dy
$$
and since $u$ is continuous in $\R^N$ and vanishes in $\Omega^c$,
there exists $z_x\in\bar \Omega$ such that $t_0=u(z_x)$, which ends
the proof.\hfill$\Box$
\medskip

%%%%%%%%%%%%PROOF%%%%%%%%%%%%%%%%%%%%%%%%%%%%%%%%%%%%%%%%%%%%%%%%%%%%%%%%%%%%%%%%%%%%%%%%%%%%%%%%%%%%%%%%%%%%%%%%%%%%%%%%%%%%%%%%%%%%%%%%%%%%%
\noindent{\bf Proof of Proposition \ref{pr 2.1}.} {\it Uniqueness}.
Let $w$ be a weak solution of
\begin{equation}\label{L0}\begin{array}{lll}
(-\Delta)^{\alpha}w=0\qquad\mbox{ in }\quad\Omega\\
\phantom{(-\Delta)^{\alpha}}w=0\qquad\mbox{ in }\quad\Omega^c.
\end{array}\end{equation}
If $\omega$ is a Borel subset of $\Omega$ and $\eta_{\omega,n}$ the
solution of
\begin{equation}\label{L00}\begin{array}{lll}
(-\Delta)^{\alpha}\eta_{\omega,n}=\zeta_n&\quad\mbox{ in }\quad\Omega\\
\phantom{(-\Delta)^{\alpha}}\eta_{\omega,n}=0&\quad\mbox{ in
}\quad\Omega^c,
\end{array}\end{equation}
where $\zeta_n:\bar\Omega\mapsto[0,1]$ is a $C^1(\bar\Omega)$
function such that
$$\zeta_n\to\chi_\omega\quad {\rm{in}}\ L^\infty(\bar\Omega)\quad {\rm{as}}\ n\to\infty.$$
Then by Lemma \ref{eta} part $(ii)$,
$\eta_{\omega,n}\in\mathbb{X}_{\alpha}$ and
$$\displaystyle\int_{\Omega}w\zeta_n dx=0.
$$
Then passing the limit of $n\to\infty$, we have
$$\displaystyle\int_{\omega}w dx=0.$$ This implies $w=0$.
\smallskip

\noindent{\it Existence and estimate (\ref{sign})}. For  $\delta>0$
we define an even convex function  $\phi_\delta$ by
\begin{equation}\label{2.3}
\phi_\delta(t)=\left\{ \arraycolsep=1pt
\begin{array}{lll}
|t|-\frac\delta2,\quad & \rm{if}\quad |t|\ge \delta ,\\[2mm]
\frac{t^2}{2\delta},\quad & \rm{if} \quad |t|< \delta/2.
\end{array}
\right.
\end{equation}
 Then for any $t,s\in\R$,  $|\phi_\delta'(t)|\le1$, $\phi_\delta(t)\to|t|$ and
$\phi_\delta'(t)\to\rm{sign}(t)$ when $\delta\to0^+$. Moreover
\begin{equation}\label{2.5}
\phi_\delta(s)-\phi_\delta(t)\ge \phi_\delta'(t) (s-t).
\end{equation}
%
%%%%%%%%%%%%%%%%%%%%%%%%%%%%%%%%%%%%%%%%%%%%%%%%%%%%%%%%%%PROOF%%%%%%%%%%%%%%%%%%%%%%%%%%%%%%%%%%%%%%%%%%%%%%%%%%%%%%%%%%%%%%%%%%%%%%%%%%%%%%%
%{\bf Proof.} Since $\phi_\delta$ is $C^{\infty}(\R)$ and$\phi_\delta(0)=0$, then $\phi_\delta(u)\in\mathbb{X}_{\alpha}$ if
%$u\in\mathbb{X}_{\alpha}$.\\ By (\ref{2.5}),  for $x\in\Omega$,
%\begin{eqnarray*}
%(-\Delta)^\alpha\phi_\delta(u)(x)&=&-\int_{\R^N}\frac{\phi_\delta(u)(x+y)-
%\phi_\delta(u)(x)}{|y|^{N+2\alpha}}dy
%\%\&\le&-\phi_\delta'(u)(x)\int_{\R^N}\frac{ u(x+y)-u
%(x)}{|y|^{N+2\alpha}}dy
%\\&=&\phi_\delta'(u)(x)(-\Delta)^\alpha  u(x).
%\end{eqnarray*}
%The proof is complete.\hfill$\Box$\\[1.5mm]
%%%%%%%%%%%%%%%%%%%%%%%%%%%%%%%%%%%%%%%%%%%%%%%%%%%%%%%%%%%%%%%%%%%%%%%%%%%%%%%%%%%%%%%%%%%%%%%%%%%%%%%%%%%%%%%%%%%%%%%%%%%%%%%%%%%%%%%%%%%%%%
Let $\{\nu_n\}$ be a sequence functions in $C^1(\bar\Omega)$ such
that
$$\lim_{n\to\infty} \int_\Omega|\nu_n-\nu|\rho^\alpha dx=0.$$
Let $u_n$ be the corresponding solution  to (\ref{homo}) with
right-hand side $\nu_n$, then by Lemma \ref{eta},
$u_n\in\mathbb{X}_{\alpha}$ and by Lemmas \ref{lm 2.1}, \ref{lm
2.2}, for any $\delta>0$ and $\xi\in\mathbb{X}_{\alpha},\ \xi\ge0$,
\begin{equation}\label{2.7}
 \arraycolsep=1pt
\begin{array}{lll}
\displaystyle\int_\Omega \phi_\delta(u_n)(-\Delta)^\alpha \xi
dx&=\displaystyle\int_\Omega
 \xi(-\Delta)^\alpha\phi_\delta(u_n) dx\\[4mm]
 &\le\displaystyle\int_\Omega
 \xi\phi_\delta'(u_n)(-\Delta)^\alpha u_n dx
\\[4mm]\displaystyle
&=\displaystyle \int_\Omega  \xi \phi_\delta'(u_n) \nu_ndx.
\end{array}
\end{equation}
Letting $\delta\to 0$, we obtain
\begin{equation}\label{L1}
\displaystyle\int_\Omega |u_n|(-\Delta)^\alpha \xi dx\leq
\displaystyle \int_\Omega  \xi {\rm {sign}} (u_n) \nu_ndx\leq
\displaystyle \int_\Omega  \xi |\nu_n|dx.
\end{equation}
If we take $\xi=\eta_1$, we derive from Lemma \ref{eta}
\begin{equation}\label{L2}
\displaystyle\int_\Omega |u_n| dx\leq  C\displaystyle \int_\Omega
|\nu_n|\rho^{\alpha}dx.
\end{equation}
Similarly
\begin{equation}\label{L3}
\displaystyle\int_\Omega |u_n-u_m| dx\leq  C\displaystyle
\int_\Omega |\nu_n-\nu_m|\rho^{\alpha}dx.
\end{equation}
Therefore $\{u_n\}$ is a Cauchy sequence in $L^1$ and its limit $u$
is a weak solution of (\ref{homo}). Letting $n\to\infty$ in
$(\ref{L1})$ we obtain (\ref{sign}).   Inequality (\ref{sign+}) is
proved by replacing $\phi_\delta$ by $\tilde\phi_\delta$ which is
zero on $(-\infty,0]$ and $\phi_\delta$ on $[0,\infty)$.
\hfill$\Box$\medskip

The next result is a higher order regularity result

\begin{proposition}\label{pr 4} Let the assumptions of Proposition \ref{general} be fulfilled and $0\leq \beta\leq\alpha$. Then for $p\in (1,\frac{N}{N+\beta-2\alpha})$ there exists $c_p>0$ such that for any $\nu\in L^{1}(\Omega,\rho^{\beta}dx)$
  \begin{equation}\label{power1'}
  \norm{\mathbb G[\nu]}_{W^{2\alpha-\gamma,p}(\Omega)}\leq c_p\norm \nu_{L^{1}(\Omega,\rho^{\beta}dx)}
\end{equation}
where $\gamma=\beta+\frac{N}{p'}$ if $\beta>0$ and
$\gamma>\frac{N}{p'}$ if $\beta=0$.
 \end{proposition}
 %%%%%%%%%%%%PROOF%%%%%%%%%%%%%%%%%%%%%%%%%%%%%%%%%%%%%%%%%%%%%%%%%%%%%%%%%%%%%%%%%%%%%%%%%%%%%%%%%%%%%%%%%%%%%%%%%%%%%%%%%%%%%%%%%%%%%%%%%%%%%
 {\bf Proof.}  We use Stampacchia's duality method \cite{St} and put $u=\mathbb G[\nu]$. If $\psi\in C^\infty_c(\bar \Omega)$, then
  \begin{equation}\label{power1-1}
\begin{array}{lll}\displaystyle\left|\int_{\Omega}\psi(-\Delta)^\alpha udx\right|\leq \int_{\Omega}|\nu||\psi| dx\\[4mm]\phantom{\left|\int_{\Omega}\psi(-\Delta)^\alpha udx\right|}\leq \displaystyle\sup_{\Omega}|\rho^{-\beta}\psi|
\int_{\Omega}|\nu|\rho^{\beta} dx\\[4mm]\phantom{\left|\int_{\Omega}\psi(-\Delta)^\alpha udx\right|}
\leq \displaystyle\norm{\psi}_{C^{\beta}(\bar \Omega)}\norm
\nu_{L^{1}(\Omega,\rho^{\beta}dx)}.
\end{array}
\end{equation}
By Sobolev-Morrey embedding type theorem (see e.g. \cite[Th
8.2]{NPV}), for any $p\in (1,\frac{N}{N+\beta-2\alpha})$ and
$p'=\frac{p}{p-1}$,
$$\norm{\psi}_{C^{\beta}(\bar \Omega)}\leq C\norm{\psi}_{W^{\gamma,p'}(\Omega)}
$$
with $\gamma=\beta+\frac{N}{p'}$ if $\beta>0$ and
$\gamma>\frac{N}{p'}$ if $\beta=0$. Therefore,
   \begin{equation}\label{power1-2}
\begin{array}{lll}\displaystyle\left|\int_{\Omega}\psi(-\Delta)^\alpha udx\right|\leq C\norm{\psi}_{W^{\gamma,p'}(\Omega)}\norm \nu_{L^{1}(\Omega,\rho^{\beta}dx)},
\end{array}
\end{equation}
which implies that the mapping $\psi\mapsto
\int_{\Omega}\psi(-\Delta)^\alpha udx$ is continuous on
$W^{\gamma,p'}(\Omega)$ and thus
   \begin{equation}\label{power1-3}\norm{(-\Delta)^\alpha u}_{W^{-\gamma,p}(\Omega)}\leq C\norm \nu_{L^{1}(\Omega,\rho^{\beta}dx)}.
\end{equation}
Since $(-\Delta )^{-\alpha}$ is an isomorphism from
$W^{-\gamma,p}(\Omega)$ into $W^{2\alpha-\gamma,p}(\Omega)$, it
follows that
   \begin{equation}\label{power1-4}\norm{u}_{W^{2\alpha-\gamma,p}(\Omega)}\leq C\norm \nu_{L^{1}(\Omega,\rho^{\beta}dx)}.
\end{equation}
\hfill$\Box$\medskip
 %Next we give a complete proof based upon Boccardo-Gallouet' ideas \cite{BG}. Since $\nu\in L^{\infty}$, $u$ is weak variational solution, that is
   %\begin{equation}\label{power2}
  %\int_{\Omega}(-\Delta)^{\frac{\alpha}{2}}u(-\Delta)^{\frac{\alpha}{2}}\xi dx= \int_{\Omega} \nu\xi dx
  %\end{equation}
 %for all $\xi\in W^{\alpha,2}(\Omega)$. For $k>0$ set $T_k(u)={\rm sgn}(u)\min\{k,|u|\}$, then $T_k\in W^{\alpha,2}(\Omega)\cap L^{\infty}(\Omega)$. Furthermore
% $$ T_k(u(x+y))-T_k(u(x))=T'_k(\tau u(x+y)+(1-\tau)u(x))u(x+y)-u(x)
 %$$
\begin{proposition}\label{pr5} Under the assumptions of Proposition \ref{pr 4} the mapping $\nu\mapsto \mathbb G[\nu]$ is compact from  $L^{1}(\Omega,\rho^\beta dx)$ into $L^{q}(\Omega)$ for any $q\in [1,\frac{N}{N+\beta-2\alpha})$.
 \end{proposition}
 %%%%%%%%%%%%PROOF%%%%%%%%%%%%%%%%%%%%%%%%%%%%%%%%%%%%%%%%%%%%%%%%%%%%%%%%%%%%%%%%%%%%%%%%%%%%%%%%%%%%%%%%%%%%%%%%%%%%%%%%%%%%%%%%%%%%%%%%%%%%%
  {\bf Proof.} By \cite[Th 6.5]{NPV} the embedding of $W^{2\alpha-\gamma,p}(\Omega)$ into $L^{q}(\Omega)$ is compact, this ends the proof.\hfill$\Box$\medskip

%%%%%%%%%%%%%%%%%%%%%%%%%%%%%%%%%%%%%%%%%%%%%%%%%%%%%%%%%%%%%%%%%%%%%%%%%%%%%%%%%%%%%%%%%%%%%%%%%%%%%%%%%%MAIN%THEOREM%%%%%%%%%%%%%%%%%%%%%%%%%%%%%%%%%%%%%%%%%%%%%%%%%%%%%%%%%%%%%%%%%%%%%%%%
%%%%%%%%%%%%%%%%%%%%%%%%%%%%%%%%%%%%%%%%%%%%%%%%%%%%%%%%%%%%%%%%%%%%%%%%%%%%%%%%%%%%%%%%%%%%%%%%%%%%%%%%%%%%%%%%%%%%%%%%%%%%%%%%%%%%%%%%%%%%%
%%%%%%%%%%%%%%%%%%%%%%%%%%%%%%%%%%%%%%%%%%%%%%
%%%%%%%%%%%%%%%%%%%%%%%%%%%%%%%%%%%%%%%%%%%%%%%%%%%%%%%%%%%%%%%%%%%%%%%%%%%%%%%%%%%%%%%%%%%%%%%%%%%%%%%%%%%%%%%%%%%%%%%%%%%%%%%%%%%%%%%%%%%%%
%%%%%%%%%%%%%%%%%%%%%%%%%%%%%%%%%%%%%%%%%%%%%%
%%%%%%%%%%%%%%%%%%%%%%%%%%%%%%%%%%%%%%%%%%%%%%%%%%%%%%%%%%%%%%%%%%%%%%%%%%%%%%%%%%%%%%%%%%%%%%%%%%%%%%%%%%%%%%%%%%%%%%%%%%%%%%%%%%%%%%%%%%%%%
%%%%%%%%%%%%%%%%%%%%%%%%%%%%%%%%%%%%%%%%%%%%%%
%%%%%%%%%%%%%%%%%%%%%%%%%%%%%%%%%%%%%%%%%%%%%%%%%%%%%%%%%%%%%%%%%%%%%%%%%%%%%%%%%%%%%%%%%%%%%%%%%%%%%%%%%%%%%%%%%%%%%%%%%%%%%%%%%%%%%%%%%%%%%

\setcounter{equation}{0}
\section{Proof of Theorem \ref{teo 1}}
Before proving the main  we give a general existence result in
$L^1(\Omega,\rho^\alpha dx)$.
%%%%%%%%%%%%%%%%%%%%%%%%%%%%%%%%%%%%%%%%%%%%%%%%%%%%%%%%%%%%%%%%%%%%%%%%%%%%%%%%%%%%%%%%%%%%%%%%%%%%%%%%%%%%%%%%%%%%%%%%%%%%%%%%%%%%%%%%%%%%%
\begin{proposition}\label{pr 2} Suppose  that  $\Omega$ is an open
bounded $C^2$ domain of $\R^N\ (N\ge2)$, $\alpha\in(0,1)$ and the
function $g:\R\to\R$ is continuous, nondecreasing and $rg(r)\ge0$
for all $r\in\R$. Then for any $f\in L^1(\Omega,\rho^\alpha dx)$
there exists a unique weak solution $u $ of (\ref{eq1.1}) with
$\nu=f$. Moreover the mapping $f\mapsto u$ is increasing.
\end{proposition}
%%%%%%%%%%%%%%%%%%%%%%%%%%%%%%%%%%%%%%%%%%%%%%%%%%%%%%%%%%%%%%%%%%%%%%%%%%%%%%%%%%%%%%%%%%%%%%
{\bf Proof.} {\it Step 1: Variational solutions. }  If $w\in L^2(\Omega)$, we denote by $\underline w$ its extension by $0$ in
$\Omega^c$ and by $W^{\alpha,2}_c(\Omega)$ the set of function in $L^2(\Omega)$ such that
$$\norm w^2_{W^{\alpha,2}_c(\Omega)}:=\int_{\R^N}|\hat{\underline w}|^2(1+|x|^\alpha)dx<\infty,
$$
where $\hat{\underline w}$ is the Fourier transform of $\underline w$.
For $\epsilon>0$ we set
$$J(w)=\frac{1}{2}\int_{\R^N}\left((-\Delta)^{\frac{\alpha}{2}}\underline w\right)^2dx+ \int_{\Omega}(j(w)+\epsilon w^2)dx,
$$
with domain $D(J)=\{w\in W_c^{\alpha,2}(\R^N) \mbox{ s.t. }j(w)\in L^1(\Omega)\}$ and $j(s)=\int_0^s g(t)dt$.
Furthermore since there holds $J(w)\geq \sigma \norm w_{W_c^{\alpha,2}}^2$ for some $\sigma>0$, the subdifferential $\partial J$ of
$J$ is a maximal monotone in the sense of Browder-Minty (see
\cite{Br} and the references therein) which satisfies $R(\partial J)=L^2(\Omega)$. Then for any $f\in L^2(\Omega)$ there exists a unique $u_\epsilon$ in the domain $D(\partial J)$ such that $\partial J(u_\epsilon)=f$.  Since for any $\psi\in W^{\alpha,2}_c(\Omega)$
$$\int_{\R^N}(-\Delta)^{\frac{\alpha}{2}}\underline w\,(-\Delta)^{\frac{\alpha}{2}}\underline \psi dx=(4\pi)^{\alpha}\int_{\R^N}\hat{\underline w}\,\hat{\underline \psi}|x|^{2\alpha}dx=\int_{\Omega}\psi (-\Delta)^{{\alpha}} \underline wdx,
$$
$$\partial J(u_\epsilon)=(-\Delta)^{\alpha}u_\epsilon+g(u_\epsilon)+2\epsilon u=f,
$$
with $u_\epsilon\in W^{2\alpha,2}_c(\Omega)$ such that $g(u_\epsilon)\in L^2(\Omega)$. This is also a consequence of \cite[Cor 2.11]{Br}.
 If $f$ is assumed
to be bounded, then $u\in C^{\alpha}(\overline\Omega)$ by \cite[Prop
1.1]{RS}. Note that more delicate variational formulations can be found in \cite{Ha1}, \cite{Ha2}. \smallskip

\noindent{\it Step 2: $L^1$ solutions. } For $n\in\N^*$ we denote by $u_{n,\epsilon}$ the solution of
\begin{equation}\label{L1-n}\begin{array}{ll}
(-\Delta)^{\alpha}u_{n,\epsilon}+g(u_{n,\epsilon})+2\epsilon u_{n,\epsilon}=f_n\qquad&\mbox{in }\,\Omega\\
\phantom{(-\Delta)^{\alpha}+g(u_{n,\epsilon})+2\epsilon u_{n,\epsilon}}u_{n,\epsilon}=0\qquad&\mbox{in }\,\Omega^c
\end{array}\end{equation}
 where $f_n={\rm sgn}(f)\min\{n,|f|\}$. By $(\ref{L1})$ with
$\xi=\eta_1$,
\begin{equation}\label{L2-n}\displaystyle\int_\Omega\left( |u_{n,\epsilon}|+(
2\epsilon|u_{n,\epsilon}|+|g(u_{n,\epsilon})|)\eta_1\right)dx\leq \displaystyle \int_\Omega
|f_n|\eta_1dx\leq\int_\Omega
|f|\eta_1dx ,
\end{equation}
and, for $\epsilon'>0$ and $m\in\N^*$,
\begin{equation}\label{L3-n}\begin{array}{ll}\displaystyle\int_\Omega\left( |u_{n,\epsilon}-u_{m,\epsilon'}|
+|g(u_{n,\epsilon})-g(u_{m,\epsilon'})|\eta_1\right)dx\\[4mm]\phantom{
+|g(u_{n,\epsilon})-g(u_{m,\epsilon'})|}
\leq \displaystyle \int_\Omega\left(
|f_n-f_m|+2\epsilon|u_{n,\epsilon}|+2\epsilon'|u_{m,\epsilon'}|\right)\eta_1dx.
\end{array}\end{equation}
Since $f_n\to f$ in $L^1(\Omega,\rho^\alpha dx)$, $\{u_{n,\epsilon}\}$ and
$\{g\circ u_{n,\epsilon}\}$ are Cauchy filters in $L^1(\Omega)$ and $L^1(\Omega,\rho^\alpha dx)$ respectively. Set
$u=\lim_{n\to\infty,\epsilon\to 0}u_{n,\epsilon}$, we derive from  the following identity valid for any
$\xi\in\mathbb X_\alpha$
$$\int_{\Omega}\left(u_{n,\epsilon}(-\Delta)^{\alpha}\xi+g(u_{n,\epsilon})\xi\right) dx=\int_{\Omega}\left(f_n-\epsilon u_{n,\epsilon}\right)\xi dx
$$
that $u$ is a solution of (\ref{eq1.1}). Uniqueness follows from
$(\ref{L1})$-$(\ref{L3-n})$, since for any $f$ and $f'$ in
$L^1(\Omega,\rho^\alpha dx)$, the any couple $(u,u')$ of weak
solutions with respective right-hand side $f$ and $f'$ satisfies
\begin{equation}\label{L5}
\displaystyle\int_\Omega\left( |u-u'|
+|g(u)-g(u')|\eta_1\right)dx\leq \displaystyle \int_\Omega
|f-f'|\eta_1dx.
\end{equation}
Finally, the monotonicity of the mapping $f\mapsto u$ follows from
$(\ref{sign+})$ thanks to which $(\ref{L5})$ is transformed into
\begin{equation}\label{L6}
\displaystyle\int_\Omega\left( (u-u')_+
+(g(u)-g(u'))_+\eta_1\right)dx\leq \displaystyle \int_\Omega
(f-f')_+\eta_1dx.
\end{equation}
\hfill$\Box$\medskip

%%%%%%%%%%%%%%%%%%%%%%%%%%%%%%%%%%%%%%%%%%%%%%%%%%%%%%%%%%%%%%%%%%%%%%%%%%%%%%%%%%%%%%%%%%%%%%%%%%%%%%%%PROOF%OF%THEOREM1%%%%%%%%%%%%%%%%%%%%%%%%%%%%%%%%%%%%%%%%%%%%%%%%%%%%%%%%%%%%%%%%%%%%%%%%%%%%%%%%%%%%%%%%%%%%%%%%%%%%%%%%%%%%%%%%%%%%%%
\noindent{\bf Proof of Theorem \ref{teo 1}.} Uniqueness follows from
$(\ref{L5})$. For existence we define
$$C_{\beta}(\bar \Omega)=\{\zeta\in C(\bar \Omega):\rho^{-\beta}\zeta\in C(\bar \Omega)\}$$
endowed with the norm
$$\norm{\zeta}_{C_{\beta}(\bar\Omega)}=\|\rho^{-\beta}\zeta\|_{C(\bar\Omega)}. $$
We consider a sequence $\{\nu_n\}\subset C^1(\bar \Omega)$ such that
$\nu_{n,\pm}\to\nu_{\pm}$ in the duality sense with $C_{\beta}(\bar
\Omega)$, which means
$$\lim_{n\to\infty}\int_{\bar \Omega}\zeta \nu_{n,\pm}dx=\int_{\bar \Omega}\zeta d\nu_{\pm}
$$
for all $\zeta\in C_{\beta}(\bar \Omega)$. It follows from the
Banach-Steinhaus theorem that $\norm{\nu_{n}}_{\mathfrak M
(\Omega,\rho^{\beta})}$ is bounded independently of $n$, therefore
\begin{equation}\label{L7}
\displaystyle\int_\Omega\left( |u_n| +|g(u_n)|\eta_1\right)dx\leq
\displaystyle \int_\Omega |\nu_n|\eta_1dx\leq C.
\end{equation}
Therefore $\norm{g(u_n)}_{\mathfrak M (\Omega,\rho^{\alpha})}$ is
bounded independently of $n$. For $\epsilon>0$, set
$\xi_{\epsilon}=(\eta_1+\epsilon)^{\frac{\beta}{\alpha}}-\epsilon^{\frac{\beta}{\alpha}}$,
which is concave in the interval $\eta(\bar\omega)$. Then, by Lemma
\ref{lm 2.2} part $(ii)$,
$$\displaystyle\begin{array}{lll}\displaystyle
(-\Delta)^{\alpha}\xi_{\epsilon}=\frac{\beta}{\alpha}(\eta_1+\epsilon)^{\frac{\beta-\alpha}{\alpha}}
(-\Delta)^{\alpha}\eta_1-\frac{\beta(\beta-\alpha)}{\alpha^2}(\eta_1+\epsilon)^{\frac{\beta-2\alpha}{\alpha}}\int_{\Omega}\frac{(\eta_1(y)-\eta_1(x))^2}{|y-x|^{N+2\alpha}}dy\\[4mm]
\phantom{(-\Delta)^{\alpha}\xi_{\epsilon}}\displaystyle \geq
\frac{\beta}{\alpha}(\eta_1+\epsilon)^{\frac{\beta-\alpha}{\alpha}},
\end{array}$$
and $\xi_\epsilon\in\mathbb{X}_{\alpha}$. Since
$$\int_{\Omega}\left(|u_n|(-\Delta)^{\alpha}\xi_{\epsilon}+|g(u_n)|\xi_{\epsilon}\right)dx\leq
\int_{\Omega}\xi_{\epsilon} d|\nu_n|,
$$
we obtain
$$\int_{\Omega}\left(|u_n|\frac{\beta}{\alpha}(\eta_1+\epsilon)^{\frac{\beta-\alpha}{\alpha}}+|g(u_n)|\xi_{\epsilon}\right)dx\leq
\int_{\Omega}\xi_{\epsilon} d|\nu_n|.
$$
If we let $\epsilon\to 0$, we obtain
$$\int_{\Omega}\left(|u_n|\frac{\beta}{\alpha}\eta_1^{\frac{\beta-\alpha}{\alpha}}+|g(u_n)|\eta_1^{\frac{\beta}{\alpha}}\right)dx\leq
\int_{\Omega}\eta_1^{\frac{\beta}{\alpha}} d|\nu_n|.
$$
By Lemma \ref{lm 2.2}, we derive the estimate
\begin{equation}\label{L9}
\int_{\Omega}\left(|u_n|\rho^{\beta-\alpha}+|g(u_n)|\rho^{\beta}\right)
dx\leq C\norm {\nu_n}_{\mathfrak M (\Omega,\rho^{\beta})}\leq C'.
\end{equation}
Since $u_n=\mathbb G[\nu_n-g(u_n)]$, it follows by $(\ref{annex
0})$, that
\begin{equation}\label{L8}
\norm {u_n}_{M^{k_{\alpha,\beta}}(\Omega,\rho^{\beta}dx)}\leq
\norm{\nu_{n}-g(u_n)}_{\mathfrak M (\Omega,\rho^{\beta})},
\end{equation}
where $k_{\alpha,\beta}$ is defined by $(\ref{power})$. By Corollary
\ref{pr5} the sequence $\{u_n\}$ is relatively compact in the
$L^q(\Omega)$ for $1\leq q<\frac{N}{N+\beta-2\alpha}$. Therefore
there exist a sub-sequence $\{u_{n_k}\}$ and some $u\in
L^1(\Omega)\cap L^q(\Omega)$ such that  $u_{n_k}\to u$ in
$L^q(\Omega)$ and almost every where in $\Omega$. Furthermore
$g(u_{n_k})\to g(u)$ almost every where. Put $\tilde
g(r)=g(|r|)-g(-|r|)$ and we note that $|g(r)|\le \tilde g(|r|)$ for
$r\in\R$ and $\tilde g$ is nondecreasing.  For $\lambda
>0$, we set $S_\lambda=\{x\in\Omega:|u_{n_k}(x)|>\lambda\}$  and
$\omega(\lambda)=\int_{S_\lambda}\rho^{\beta}dx$. Then for any Borel
set $E\subset\Omega$, we have
$$\displaystyle\begin{array}{lll}
\displaystyle\int_{E}|g(u_{n_k})|\rho^{\beta}dx=\int_{E\cap
S^c_{\lambda}}|g(u_{n_k})|\rho^{\beta}dx+\int_{E\cap S_{\lambda}}
|g(u_{n_k})|\rho^{\beta}dx\\[4mm]\phantom{\int_{E}|g(u_{n_k})|\rho^{\beta}dx}
\displaystyle\leq \tilde g(\lambda)\int_E\rho^{\beta}dx+\int_{S_{\lambda}}\tilde g(|u_{n_k}|)\rho^{\beta}dx\\[4mm]\phantom{\int_{E}|g(u_{n_k})|\rho^{\beta}dx}
\displaystyle\leq \tilde
g(\lambda)\int_E\rho^{\beta}dx-\int_{\lambda}^\infty \tilde
g(s)d\omega(s).
\end{array}$$
But
$$\int_{\lambda}^\infty \tilde g(s)d\omega(s)=\lim_{T\to\infty}\int_{\lambda}^T \tilde g(s)d\omega(s).
$$
Since $u_{n_k}\in M^{k_{\alpha,\beta}}(\Omega,\rho^{\beta}dx)$,
$\omega(s)\leq cs^{-k_{\alpha,\beta}}$ and
$$\displaystyle\begin{array}{lll}
\displaystyle-\int_{\lambda}^T \tilde g(s)d\omega(s) =-\left[\tilde
g(s)\omega(s)\!\!\!\!\!\!\!\!^{\phantom{\frac{X^X}{X}}}\right]_{s=\lambda}^{s=T}+\int_{\lambda}^T
\omega(s)d\tilde g(s)
\\[4mm]\phantom{\int_{\lambda}^T \tilde g(s)d\omega(s)}\displaystyle
\leq\tilde g(\lambda)\omega(\lambda)-\tilde
g(T)\omega(T)+c\int_{\lambda}^T s^{-k_{\alpha,\beta}}d\tilde g(s)
\\[4mm]\phantom{\int_{\lambda}^T \tilde g(s)d\omega(s)}\displaystyle
\leq \tilde g(\lambda)\omega(\lambda)-\tilde g(T)\omega(T)+
c\left(T^{-k_{\alpha,\beta}}\tilde
g(T)-\lambda^{-k_{\alpha,\beta}}\tilde g(\lambda)\right)
\\[4mm]\phantom{-----------\int_{\lambda}^T \tilde g(s)d\omega(s)}\displaystyle
+\frac{c}{k_{\alpha,\beta}+1}\int_{\lambda}^T
s^{-1-k_{\alpha,\beta}}\tilde g(s)ds.
\end{array}$$
By assumption $(\ref{1.4})$  there exists $\{T_n\}\to\infty$ such
that $T_n^{-k_{\alpha,\beta}}\tilde g(T_n)\to 0$ when $n\to\infty$.
Furthermore $\tilde g(\lambda)\omega(\lambda)\leq
c\lambda^{-k_{\alpha,\beta}}\tilde g(\lambda)$, therefore
$$-\int_{\lambda}^\infty \tilde g(s)d\omega(s)\leq \frac{c}{k_{\alpha,\beta}+1}\int_{\lambda}^\infty s^{-1-k_{\alpha,\beta}}\tilde g(s)ds.
$$
Notice that the above quantity on the right-hand side tends to $0$
when $\lambda\to\infty$. The conclusion follows: for any
$\epsilon>0$ there exists $\lambda>0$ such that
$$\frac{c}{k_{\alpha,\beta}+1}\int_{\lambda}^\infty s^{-1-k_{\alpha,\beta}}\tilde g(s)ds\leq \frac{\epsilon}{2}
$$
and $\delta>0$ such that
$$\int_E\rho^{\beta}dx\leq \delta\Longrightarrow \tilde g(\lambda)\int_E\rho^{\beta}dx\leq\frac{\epsilon}{2}.
$$
This proves that $\{g\circ u_{n_k}\}$ is uniformly integrable in
$L^1(\Omega,\rho^\beta dx)$. Then $g\circ u_{n_k}\to g\circ u$ in
$L^1(\Omega,\rho^\beta dx)$ by Vitali convergence theorem. Letting
$n_k\to\infty$ in the identity
$$\int_{\Omega}\left(u_{n_k}(-\Delta)^{\alpha}\xi+\xi g\circ u_{n_k}\right) dx=\int_{\Omega}\nu_{n_k}\xi dx
$$
where $\xi\in \mathbb X_{\alpha}$, it infers that $u$ is a weak
solution of $(\ref{eq1.1})$.

The right-hand side of estimate $(\ref{1.4})$ follows from the fact
that $v_{n,+}:=\mathbb G[\nu_{n,+}]$ satisfies
$$(-\Delta)^{\alpha}v_{n,+}+g(v_{n,+})=\nu_{n,+}+g(v_{n,+})\geq \nu_{n}
$$
Therefore $v_{n,+}\geq u_n$ by Proposition \ref{pr 2}. Letting
$n\to\infty$ yields to $(\ref{1.5})$. The left-hand side is proved
similarly.

% \noindent{\it   To prove the uniqueness.}  Let $u_1,u_2$ are two solutions of
%(\ref{eq1.1}) and $w=u_1-u_2$. Then $(-\Delta)^\alpha
%w=g(u_2)-g(u_1),$ where $g(u_2)-g(u_1)\in L^1(\Omega,\rho^\alpha
%dx)$. By (\ref{sign}),  for $\xi\in\mathbb{X}$, $\xi\ge0$,
%\begin{eqnarray*}
%\int_\Omega |w|(-\Delta)^\alpha \xi
%dx+\int_\Omega[g(u_1)-g(u_2)]\rm{sign}(w)\xi dx\le0.
%\end{eqnarray*}
%Together with $\int_\Omega[g(u_1)-g(u_2)]\rm{sign}(w)\xi dx\ge0$,
%then we have
%$$w=0,\quad \rm{a.e.\ in}\ \Omega.$$

\noindent{\it  To prove the mapping $\nu\mapsto u$ is increasing.}
Let $\nu_1,\nu_2\in\mathfrak{M}(\Omega,\rho^\beta)$ and $\nu_1\ge
\nu_2$, then there exist two sequences $\{\nu_{1,n}\}$ and
$\{\nu_{2,n}\}$ in $C^\infty(\bar\Omega)$ such that $\nu_{1,n}\ge
\nu_{2,n}$ and
$$\nu_{i,n}\to\nu_i\quad \rm{as}\ n\to\infty,\quad i=1,2.$$ Let
$u_{i,n}$ be the unique solution of (\ref{eq1.1}) with  $\nu_{i,n}$
and $u_{i}$ be the unique solution of (\ref{eq1.1}) with  $\nu_{i}$
where $i=1,2$. Then $u_{1,n}\ge u_{2,n}$. Moveover,  by uniqueness
 $u_{i,n}$ convergence to $u_i$ in $L^1(\Omega)$ for $i=1$
and $i=2$. Then we have $u_1\ge u_2$. \hfill$\Box$\medskip

%%%%%%%%%%%%%%%%%%%%%%%%%%%%%%%%%%%%%%%%%%%%%%%%%%%%%%%%%%%%%%%%%%%%%%%%%%%%%%%%%%%%%%%%%%%%%%%%%%%%%%%%%%%%%%%%%%%%%%%%%%%%%%%%%%%%%%%%%%%%%%%STABILITY%%%%%%%%%%%%%%%%%%%%%%%%%%%%%%%%%%%%%%%%%%%%%%%%%%%%%%%%%%%%%%%%%%%%%%%%%%%%%%%%%%%%%%%%%%%%%%%%%%%%%%%%%%%%%%%%%%%%%%%%%%%%%%%%%%%%%%%%%%%%%%%%%%%%%%%%%%%%%%%%%%%%%

\begin{corollary}\label{stability} Under the hypotheses of Theorem \ref{teo 1}, we further
assume that  $\{\nu_n\}$ is a sequence of measures in $\mathfrak M
(\Omega,\rho^{\beta})$ and $\nu\in \mathfrak M
(\Omega,\rho^{\beta})$ such that  for any $\xi\in C_\beta(\bar
\Omega)$,
$$\int_{\Omega}{}\xi d\nu_n\to \int_{\Omega}{}\xi d\nu\quad{\rm{as}}\quad n\to\infty.
$$
Then the sequence $\{u_{n}\}$ of weak solutions to
\begin{equation}\label{eq1.1n}
 \arraycolsep=1pt
\begin{array}{lll}
 (-\Delta)^\alpha  u_{n}+g\circ u_{n}=\nu_n\quad & \rm{in}\quad\Omega,\\[2mm]
 \phantom{ (-\Delta)^\alpha  +g\circ u_{n}}
u_{n}=0\quad & \rm{in}\quad \Omega^c,
\end{array}
\end{equation}
converges to the solution $u$ of $(\ref{eq1.1})$ in $L^q(\Omega)$
for $1\leq q<\frac{N}{N+\beta-2\alpha}$ and $\{g\circ u_{n}\}$
converges to $g\circ u$ in $L^1(\Omega,\rho^{\beta}dx)$.
\end{corollary}
%%%%%%%%%%%%%%%%%%%%%%%%%%%%%%%%%%%%%%%%%%%%%%%%%%%%%%%%%%%%%%%%%%%%%%%%%%%%%%%%%%%%%%%%%%%%%%
{\bf Proof.} The method is an adaptation of \cite{V1}. Since
$\nu_n\to\nu$ in the duality sense of $C_\beta(\overline \Omega)$,
there exists $M>0$ such that
$$\norm{\nu_n}_{\mathfrak M(\Omega,\rho^{\beta})}\leq M,\qquad\forall n\in\mathbb N.
$$
Therefore $(\ref{L9})$ and $(\ref{L8})$ hold (but with $u_n$
solution of $(\ref{eq1.1n})$). The above proof shows that $\{g\circ
u_n\}$ is uniformly integrable in $L^1(\Omega,\rho^{\beta}dx)$ and
$\{ u_n\}$ relatively compact in $L^q(\Omega)$ for $1\leq
q<\frac{N}{N+\beta-2\alpha}$. Thus, up to a subsequence $\{
u_{n_k}\}\subset\{ u_n\}$, $u_{n_k}\to u$, and $u$ is the weak
solution of $(\ref{eq1.1})$. Since $u$ is unique, $u_{n}\to u$ as
$n\to\infty$. \hfill$\Box$\medskip
%%%%%%%%%%%%%%%%%%%%%%%%%%%%%%%%%%%%%%%%%%%%%%%%%%%%%%%%%%%%%%%%%%%%%%%%%%%%%%%%%%%%%%%%%%%%%%%%%%%%%%%%%%%%%%%OLD%%%%%%%%%%%%%%%%%%%%%%%%%%%%%%%%%%PROOF%%%%%%%%%%%%%%%%%%%%%%%%%%%%%%%%%%%%%%%%%%%%%%%%%%%%%%%%%%%%%%%%%%%%%%%%%%%%%%%%%%%%%%%%%%%%%%%%%%%%%%%%%%%%%%%%%%%%%%%%%%%%%%%%%%%%%%%%%%%%%%%%%%%%%%%%%%%%%%%%%%%%
\begin{remark}
Under the hypotheses of Theorem \ref{teo 1}, we assume $\nu\ge0$,
then
\begin{equation}\label{4.5}
\mathbb{G}(\nu)-\mathbb{G}(g(\mathbb{G}(\nu)))\le u\le
\mathbb{G}(\nu).
\end{equation}
Indeed, since $g$ is nondecreasing and $u\le \mathbb{G}(\nu)$, then
\begin{eqnarray*}
u&=& \mathbb{G}(\nu)- \mathbb{G}(g(u))
\\&\ge& \mathbb{G}(\nu)- \mathbb{G}(g(\mathbb{G}(\nu))).
\end{eqnarray*}
\end{remark}

%%%%%%%%%%%%%%%%%%%%%%%%%%%%%%%%%%%%%%%%%%%%%%%%%%%%%%%%%%%%%%%%%%%%%%%%%%%%%%%%%%%%%%%%%%%%%%%%%%%%%%%%%%%%%%%%%%%%%%%%%%
%%%%%%%%%%%%%%%%%%%%%%%%%%%%%%%%%%%%%%%%%%%%%%%%%%%%%%%%%%%%%%%%%%%%%%%%%%%%%%%%%%%%%%%%%%%%%%%%%%%%%%%%%%%%%%%%%%%%%%%%%%
\setcounter{equation}{0}
\section{Applications}

\subsection{The case of a Dirac mass}

In this subsection we characterize the asymptotic behavior of a
solution near a singularity created by a Dirac mass.

\begin{teo}\label{teo 2}
Assume that  $\Omega$ is an open,  bounded and $C^2$ domain of
$\R^N\ (N\ge2)$ with $0\in\Omega$, $\alpha\in(0,1)$, $\nu=\delta_0$
and the function $g:[0,\infty)\to[0,\infty)$ is  continuous,
nondecreasing and (\ref{1.4}) holds for
\begin{equation}\label{4.1}
k_{\alpha,0}=\frac N{N-2\alpha}.
\end{equation}
Then problem (\ref{eq1.1}) admits a unique positive weak solution
$u$ such that
\begin{equation}\label{4.30}
\lim_{x\to0}u(x)|x|^{N-2\alpha}=C,
 \end{equation}
for some $C>0$.
 \end{teo}

 \begin{remark}
We note here that a weak solution $u$ of (\ref{eq1.1}) with
$\nu=\delta_0$ satisfies
\begin{equation}\label{eq4.1}
 \arraycolsep=1pt
\begin{array}{lll}
 (-\Delta)^\alpha  u+g(u)=0\quad & \rm{in}\quad\Omega\setminus\{0\},\\[2mm]
 \phantom{ (-\Delta)^\alpha  +g(u)}
u=0\quad & \rm{in}\quad \R^N\setminus\Omega.
\end{array}
\end{equation}
The asymptotic behavior $(\ref{4.30})$ is one of the possible
singular behaviors of solutions  of $(\ref{eq4.1})$ given in
\cite{CV1}.
\end{remark}

Before proving Theorem \ref{teo 2}, we give an auxiliary lemma.
\begin{lemma}\label{lm 4.1}
Assume that $g:[0,\infty)\to[0,\infty)$ is  continuous,
nondecreasing and (\ref{1.4}) holds with $k_{\alpha,\beta}>1$. Then
$$\lim_{s\to\infty}g(s)s^{-k_{\alpha,\beta}}=0.$$

\end{lemma}
{\bf Proof.} Since
\begin{eqnarray*} \int_s^{2s}g(t)t^{-1-k_{\alpha,\beta}}dt\ge
g(s)(2s)^{-1-k_{\alpha,\beta}}\int_s^{2s}dt=2^{-1-k_{\alpha,\beta}}g(s)s^{-k_{\alpha,\beta}}
\end{eqnarray*}
and by (\ref{1.4}),
\begin{eqnarray*}
\lim_{s\to\infty}\int_s^{2s}g(t)t^{-1-k_{\alpha,\beta}}dt=0.
\end{eqnarray*}
Then
$$\lim_{s\to\infty}g(s)s^{-k_{\alpha,\beta}}=0.$$
The proof is complete.\hfill$\Box$\\[2mm]

%%%%%%%%%%%%%%%%%%%%%%%%%%%%%%%%%%%%%%%%%%%%%%%%%%%%%%%%%%%%%%%%%%%%%%%%%%%%%%%%%%%%%%%%%%%%%%%%%%%%%%%%%%%%%%%%%%%%%%%%%%%%%%%%%%%%%%%%%%%%%%%%%%%%%%%%%%%%%%%%%%%%%%%%%%%%%%%%%%%%%%%%%%%%
\noindent{\bf Proof of Theorem \ref{teo 2}.} Existence, uniqueness
and positiveness follow from Theorem \ref{teo 1} with $\beta=0$. For
(\ref{4.30}), we shall use (\ref{1.5}).  From \cite{CV0} there
holds,
\begin{equation}\label{4.3}
0<\frac{C}{|x|^{N-2\alpha}}-G(x,0)<\frac{C}{\rho(0)^{N-2\alpha}},\quad
x\in\Omega\setminus\{0\}.
\end{equation}
 for some $C>0$ dependent of $N$ and $\alpha$. Since
$$\mathbb{G}(\delta_0)(x)=G(x,0)<\frac{C}{|x|^{N-2\alpha}},\quad x\in\Omega\setminus\{0\},$$
then
\begin{eqnarray*}
0&\le&\mathbb{G}(g(\mathbb{G}(\delta_0)))(x)|x|^{N-2\alpha}
\\&\le&\int_\Omega
\frac1{|x-y|^{N-2\alpha}}g(\frac{C}{|y|^{N-2\alpha}})dy|x|^{N-2\alpha}
\\&\le&\int_\Omega
\frac1{|e_x-y|^{N-2\alpha}}g(\frac{C}{(|x||z|)^{N-2\alpha}})dz|x|^N
\\&=&|x|^N\int_{\Omega\cap B_{1/2}(e_x)}
\frac1{|e_x-y|^{N-2\alpha}}g(\frac{C}{(|x||z|)^{N-2\alpha}})dz\\&&+|x|^N\int_{\Omega\cap
B_{1/2}^c(e_x)}
\frac1{|e_x-y|^{N-2\alpha}}g(\frac{C}{(|x||z|)^{N-2\alpha}})dz
\\&:=&A_1(x)+A_2(x),\quad x\in\Omega\setminus\{0\},
\end{eqnarray*}
 where $e_x=x/|x|$. By Lemma \ref{lm 4.1},
\begin{eqnarray*}
A_1(x)&\le&|x|^Ng(\frac{2^{N-2\alpha}C}{|x|^{N-2\alpha}})\int_{
B_{1/2}(e_x)} \frac1{|e_x-y|^{N-2\alpha}}dz
\\&\to&0\quad{\rm{as}}\ \  |x|\to0,
\end{eqnarray*}
and by (\ref{1.4}),
\begin{eqnarray*}
A_2(x)&\le&\bar C|x|^N\int_{B_R(0)}
g(\frac{C}{(|x||z|)^{N-2\alpha}}) dz
\\&\le&\bar C\int_{\frac{R^{1/(N-2\alpha)}}{|x|}}^\infty g(Cs)s^{-1-\frac
N{N-2\alpha}}ds
\\&\to&0\quad{\rm{as}}\ \ |x|\to0,
\end{eqnarray*}
where $R>0$ such that $B_R(0)\supset \Omega$. That is
\begin{equation}\label{4.4}
\lim_{|x|\to0}\mathbb{G}(g(\mathbb{G}(\delta_0)))(x)|x|^{N-2\alpha}=0.
\end{equation}
We plug (\ref{4.3}) and (\ref{4.4}) into (\ref{4.5}), then
(\ref{4.30}) holds. \hfill$\Box$

\subsection{The power case}
If $g(s)=|s|^{k-1}s$ with $k\geq 1$, then $(\ref{1.4})$ is satisfied
if $1\leq k<k_{\alpha,\beta}$ where $k_{\alpha,\beta}$ defined by
$(\ref{power})$ is called the {\it critical exponent} with limit
values $ k_{\alpha,0}=\frac{N}{N-2\alpha}$ and
$k_{\alpha,\alpha}=\frac{N+\alpha}{N-\alpha}$. If we consider the
problem
\begin{equation}\label{P1}\begin{array}{ll}
(-\Delta )^{\alpha}u+|u|^{k-1}u=\nu\qquad\mbox{in }\quad \Omega,\\[2mm]\phantom{(-\Delta )^{\alpha}u+|u|^{q-1}}
u=0\qquad\mbox{in }\quad\Omega^c,
\end{array}\end{equation}
then if $1<k<k_{\alpha,\beta}$ it is solvable for any $\nu\in
\mathfrak M(\Omega,\rho^{\beta})$, but it may not be the case if
$k\geq k_{\alpha,\beta}$. As in the case $\alpha=1$, the sharp
solvability of $(\ref{P1})$ is associated to a concentration
property of the measure $\nu$ and this concentration is expressed by
the mean of Bessel capacities. If $k>1$ and $k'=\frac{k}{k-1}$, we
define for any compact set $K\subset \Omega$,
\begin{equation}\label{P2}\begin{array}{ll}
C^\Omega_{2\alpha,k'}(K)=\inf\{\norm{\phi}^{k'}_{W^{2\alpha,k'}(\Omega)}:\phi\in
C^{\infty}_c(\Omega),0\leq \phi\leq 1,\phi\equiv 1\mbox{ on }K\}.
\end{array}\end{equation}
Then $C_{2\alpha,k'}$ is an outer measure or capacity in $\Omega$
extended to Borel sets by standard processes. Our result is the
following in the case of bounded measures
%%%%%%%%%%%%%%%%%%%%%%%%%%%%%%%%%%%%%%%%%%%%%%%%%%%%%%
%%%%%%%%%%%%%%%%%%%THEOREM%%%%%%%%%%%%%%%%%%%%%%%%%%%%%%%%%%%%%%%%%%%%%%%%%%%%%%%%%%%%%%%%%%%%%%%%%%%%%%%%%%%%%%%%%%%

\begin{teo}\label{crit1} Assume $\Omega$ is an open bounded $C^2$ domain in $\R^N$ and $k>1$. Then problem $(\ref{P1})$ can be solved with a nonnegative bounded measure $\nu$ if and only if $\nu$ satisfies on compact subsets $K\subset\Omega$
\begin{equation}\label{P3}\begin{array}{ll}
C^\Omega_{2\alpha,k'}(K)=0\Longrightarrow \nu(K)=0.
\end{array}\end{equation}
\end{teo}
{\bf Proof.} {\it 1-The condition is necessary}. Assume $u$ is a
weak solution and let $K\subset\Omega$ be compact. Let $\phi\in
C^{\infty}_c(\R^N)$ such that $0\leq \phi\leq 1$ and $\phi(x)=1$
for all $x\in K$, and set $\xi=\phi^{k'}$, then $\xi\in\mathbb
X_\alpha$ and
$$\int_{\Omega}\left(u(-\Delta)^{\alpha}\xi+u^k\xi\right)dx=\int_{\Omega}\xi d\nu.
$$
Since $\xi\geq\chi_{K}$ it follows from $(\ref{K1})$ that
\begin{equation}\label{P3'}\int_{\Omega}\left(k'\phi^{k'-1}u(-\Delta)^{\alpha}\phi+\phi^{k'}u^{k}\right)dx\geq \nu(K).
\end{equation}
By H\"older's inequality
\begin{equation}\label{P4}\left |\int_{\Omega}\phi^{k'-1}u(-\Delta)^{\alpha}\phi dx\right |
\leq  \left(\int_{\Omega}\phi^{k'}u^{k}
dx\right)^{\frac{1}{k}}\left(\int_{\Omega}\left|(-\Delta)^{\alpha}\phi\right|^{k'}
dx\right)^{\frac{1}{k'}}
\end{equation}
By \cite[Th 5.4]{NPV}, there exists $\tilde \phi\in W^{2\alpha,k'}(\R^N)$ such that $\tilde\phi\lfloor_{\Omega}=\phi$ and
$$\norm{\tilde\phi}_{W^{2\alpha,k'}(\R^N)}\leq C\norm{\phi}_{W^{2\alpha,k'}(\Omega)}
$$
Then, by standard regularity result on the Riesz potential $(-\Delta)^{-{\alpha}}$ in $\R^N$,
\begin{equation}\label{P5}\begin{array}{lll}\displaystyle
\left |\int_{\Omega}\phi^{k'-1}u(-\Delta)^{\alpha}\phi dx\right |
\leq  \left(\int_{\Omega}\phi^{k'}u^{k}
dx\right)^{\frac{1}{k}}\left(\int_{\R^N}\left|(-\Delta)^{\alpha}\phi\right|^{k'}dx\right)^{\frac{1}{k'}}\\[4mm]
\phantom{\left |\int_{\Omega}\phi^{k'-1}u(-\Delta)^{\alpha}\phi dx\right |}
\displaystyle\leq C'\left(\int_{\Omega}\phi^{k'}u^{k}dx\right)^{\frac{1}{k'}}\norm{\tilde\phi}_{W^{2\alpha,k'}(\R^N)}\\[4mm]
\phantom{\left |\int_{\Omega}\phi^{k'-1}u(-\Delta)^{\alpha}\phi dx\right |}
\displaystyle\leq C'\left(\int_{\Omega}\phi^{k'}u^{k}dx\right)^{\frac{1}{k'}}\norm{\phi}_{W^{2\alpha,k'}(\Omega)}.
\end{array}\end{equation}
Therefore $(\ref{P5})$, yields to
\begin{equation}\label{P6}
C\norm{\phi}_{W^{2\alpha,k'}(\Omega)}\left(\int_{\Omega}\phi^{k'}u^{k}dx\right)^{\frac{1}{k}}+\int_{\Omega}\phi^{k'}u^{k}dx\geq
\nu(K).
\end{equation}
If $C^\Omega_{2\alpha,k'}(K)=0$, there exists a sequence
$\{\phi_n\}\subset C^{\infty}_c(\Omega)$ such that $0\leq \phi_n\leq
1$ and $\phi_n=1$ on $K$ and
$\norm{\phi_n}_{W^{2\alpha,k'}(\Omega)}\to 0$ as $n\to\infty$.
Furthermore $K$ has zero Lebesgue measure and $\phi_n\to 0$ almost
everywhere. If we replace $\phi$ by $\phi_n$ in $(\ref{P6})$ and let
$n\to\infty$ we obtain $\nu(K)=0$.\smallskip

\noindent {\it 2-The condition is sufficient}. We first assume that
$\nu\in W^{-2\alpha,k}(\Omega)\cap \mathfrak M^b_+(\Omega)$;  for
$n\in\N$, we denote by $u_n$ the solution of
\begin{equation}\label{Pn}\begin{array}{ll}
(-\Delta )^{\alpha}u+|T_n(u)|^{k-1}T_n(u)=\nu\qquad\mbox{in }\quad\Omega\\[2mm]
\phantom{(-\Delta )^{\alpha}+|T_n(u)|^{k-1}T_n(u)} u=0\qquad\mbox{in
}\quad\Omega^c
\end{array}\end{equation}
where $T_n(r)={\rm sign}(r)\min\{n,|r|\}$. Such a solution exists by
Theorem \ref{teo 1}, is nonnegative and the sequence $\{u_n\}$ is
decreasing and converges to some nonnegative $u$ since $\{T_n(r)\}$
is increasing on $\R_+$. Furthermore
$$0\leq u_n\leq\mathbb G[\nu],
$$
by $(\ref{1.5})$. This implies that the convergence holds in
$L^1(\Omega)$. Since $\nu\in W^{-2\alpha,k}(\Omega)$, $G[\nu]\in
L^k(\Omega)$, it infers that
$$|T_n(u_n)|^{k-1}T_n(u_n)=(T_n(u_n))^k\leq (\mathbb G[\nu])^k.
$$
Since for any $\xi\in\mathbb X_\alpha$ there holds
\begin{equation}\label{P7}
\int_\Omega\left(u_n(-\Delta)^{\alpha}\xi+(T_n(u_n))^k\xi\right)dx=\int_\Omega\xi
d\nu
\end{equation}
we can let $n\to\infty$ and conclude that $u$ is a solution of
$(\ref{P1})$, unique by $(\ref{L5})$. Next we assume that
$(\ref{P3})$ holds. By a result of Feyel and de la Pradelle
\cite{FdP} (see also \cite{DM}), there exists an increasing sequence
$\{\nu_n\}\subset W^{-2\alpha,k}(\Omega)\cap \mathfrak
M_+^b(\Omega)$ which converges to $\nu$ in the weak sense of
measures.  This implies that the sequence $\{u_n\}$ of weak
solutions of
\begin{equation}\label{P1n}\begin{array}{ll}
(-\Delta )^{\alpha}u_n+u_n^{k}=\nu_n\qquad&\mbox{in }\quad\Omega\\[2mm]
\phantom{(-\Delta )^{\alpha}+u_n^{k}} u_n=0\qquad&\mbox{in
}\quad\Omega^c
\end{array}\end{equation}
is increasing with limit $u$. Taking $\eta_1:=\mathbb G[1]$ as a
test function in the weak formulation, we have
$$\int_\Omega\left(u_n+u_n^{k}\eta_1\right)dx=\int_\Omega\eta_1 d\nu_n\leq \int_\Omega\eta_1 d\nu.
$$
Therefore $u_n\to u$ in $L^1(\Omega)\cap
L^k(\Omega,\rho^{\alpha}dx)$. Letting $n\to \infty$ we deduce that
$u$ satisfies $(\ref{P1})$.\hfill$\Box$\medskip

\begin{remark}\label{gene} If $\nu$ is a signed bounded measure a sufficient condition for solving $(\ref{P1})$
is
\begin{equation}\label{P'3}\begin{array}{ll}
C^\Omega_{2\alpha,k'}(K)=0\Longrightarrow |\nu|(K)=0.
\end{array}\end{equation}
This can be obtained by using the fact that the solutions of
$(\ref{P1})$ with right-hand side $\nu_+$ and $-\nu_-$  are
respectively a supersolution and a subsolution of $(\ref{P1})$. It
is not clear whether it is also a necessary condition.
\end{remark}
%%%%%%%%%%%%%%%%%%%%%%%%%%%%%%%%%%%%%%%%%%%%%%%%%%%%%%%%%%%%%%%%%%%%%%%%%%%%%%%%%%%%%%%%%%%%%%%%%%%%%%%%%%%%%%%%%%%%%
%%%%%%%%%%%%%%%%%%%%%%%%%%%%%%%%%%%%%%%%%%%%%%%%%%%%%%%%%%%%%%%%%%%%%%%%%%%%%%%%%%%%%%%%%%%%%%%%%%%%%%%%%%%%%%%%%%%%%


\begin{thebibliography}{99}

\bibitem {BP2} P. Baras and M. Pierre, Singularit\'{e} s\'{e}liminables pour des \'{e}quations
semi lin\'{e}aires, {\it Ann. Inst. Fourier Grenoble 34} , 185-206
(1984).

\bibitem {BB11} Ph. B\'{e}nilan and H. Brezis, Nonlinear problems related to the
Thomas-Fermi equation, {\it J. Evolution Eq. 3}, 673-770, (2003).

\bibitem {BBC} Ph. B\'{e}nilan, H. Brezis and M. Crandall, A semilinear elliptic equation in $L^1(\R^N )$, {\it Ann. Sc. Norm. Sup. Pisa Cl. Sci. 2},523-555 (1975).

\bibitem {Hung} M. F. Bidaut-V\'{e}ron,  N. Hung and
L. V\'{e}ron, Quasilinear Lane-Emden equations with absorption and
measure data, {\it arXiv:1212.6314v2 [math.AP],} 15 (Jan 2013).

\bibitem {BV} M. F. Bidaut-V\'{e}ron  and L. Vivier, An elliptic semilinear equation
with source term involving boundary measures: the subcritical case,
{\it Rev. Mat. Iberoamericana 16}, 477-513 (2000).

\bibitem {Br} H. Brezis, Op\'erateurs maximaux monotones et semi-groupes de contractions dans les espaces de Hilbert, {\it North-Holland Mathematics Studies. 5. Notas de matematica 50}, North-Holland, Amsterdam (1973).

\bibitem {B12} H. Brezis, Some variational problems of the Thomas-Fermi type.
Variational inequalities and complementarity problems, {\it Proc.
Internat. School, Erice, Wiley, Chichester}, 53-73 (1980).


\bibitem {CT}
X. Cabr\'{e} and Y. Sire, Nonlinear equations for fractional Laplacians I: Regularity, maximum principles, and Hamiltonian estimates, {\it Ann. I. H. Poincar\'{e} – AN, } available online 7 February (2013).

\bibitem {CS0}
L. Caffarelli and L. Silvestre, An extension problem related to the
fractional laplacian, {\it Comm. Partial Differential Equations 32},
1245-1260 (2007).

\bibitem {CS1}
 L. Caffarelli and  L. Silvestre, Regularity theory for fully non-linear integrodifferential equations,
 {\it Communications on Pure and Applied Mathematics 62}, 597-638.(2009)


\bibitem {CS} Z. Chen, and R. Song, Estimates on Green functions and poisson kernels for symmetric stable process, {\it Math.
Ann. 312}, 465-501 (1998).

\bibitem {CV0} H. Chen and L. V\'{e}ron, Solutions of fractional equations involving sources and Radon
measures, {\it preprint.}

\bibitem {CV1} H. Chen and L. V\'{e}ron, Singular solutions of fractional elliptic equations with
absorption, {\it arXiv:1302.1427v1, [math.AP]}, 6 (Feb 2013).


\bibitem {CFQ} H. Chen, P. Felmer and A. Quaas, Large solution to elliptic  equations involving fractional Laplacian,
{\it preprint.}

\bibitem {CLO2}
W. Chen, C. Li and B. Ou, Classification of solutions for an
integral equation, {\it Comm. Pure Appl. Math. 59}, 330-343 (2006).

\bibitem {CC} R. Cignoli  and M. Cottlar, An Introduction to Functional Analysis,
{\it North-Holland, Amsterdam}, 1974.


\bibitem {DM} G. Dal Maso, On the integral representation of certain local functionals, {\it Ricerche Mat. 32}, 85-113 (1983).

\bibitem {FT}  P. Felmer and E. Topp, Uniform equicontinuity for a family of zero order operators approaching the fractional Laplacian, {\it preprint.}

\bibitem {FdP} D. Feyel and A. de la Pradelle, Topologies fines et compactifications associ\'{e}es \`{a} certains espaces de
Dirichlet,{\it  Ann. Inst. Fourier (Grenoble) 27}, 121-146 (1977).

\bibitem {GV} A. Gmira  and L. V\'{e}ron , Boundary singularities of solutions of some
nonlinear elliptic equations, {\it Duke Math. J.  64}, 271-324
(1991).

\bibitem {Ha1} H. Hajeij,  Existence of minimizers of functionals involving the fractional gradient in the absence of compactness, 
{ J. Math. Anal. Appl. 399}, 17Ð26 (2013).

\bibitem {Ha2} H. Hajeij Variational problems related to some fractional kinetic equations, arXiv:1205.1202, preprint.

\bibitem {KPU} K. Karisen, F. Petitta  and S. Ulusoy, A duality approach to the fractional laplacian with measure data,
{\it Publ. Mat. 55}, 151-161 (2011).

\bibitem {L} Y.Y. Li, Remark on some conformally invariant integral equations:
the method for moving spheres, {\it J. Eur. Math. Soc. 6}, 153-180
(2004).

\bibitem {MV1} M. Marcus  and L. V\'{e}ron, The boundary trace of positive solutions of
semilinear elliptic equations: the subcritical case, {\it Arch. Rat.
Mech. Anal. 144},  201-231 (1998).


\bibitem {MV2} M. Marcus  and L. V\'{e}ron, The boundary trace of positive
solutions of semilinear elliptic equations: the supercritical case,
{\it J. Math. Pures Appl. 77}, 481-524 (1998).

\bibitem {MV3} M. Marcus  and L. V\'{e}ron, Removable singularities and boundary
traces, {\it J. Math. Pures Appl. 80}, 879-900 (2001).

\bibitem {MV4} M. Marcus  and L. V\'{e}ron, The boundary trace and generalized B.V.P. for
semilinear elliptic equations with coercive absorption, {\it Comm.
Pure Appl. Math. 56}, 689-731(2003).


\bibitem{NPV}  E. Di Nezza, G. Palatucci and E. Valdinoci,  Hitchhiker's guide to the fractional Sobolev spaces,
  {\it arXiv:1104.4345v3, [math.FA]}, 19 (Nov 2011).

\bibitem {RS} X. Ros-oton and J. Serra, The Dirichlet problem for the fractional
laplacian: regularity up to the boundary, {\it arXiv:1207.5985v1
[math.AP]}, 25 (Jul 2012).


\bibitem {S}
L. Silvestre, Regularity of the obstacle problem for a fractional
power of the laplace operator, {\it Comm. Pure Appl. Math. 60},
67-112 (2007).

\bibitem {S1}
Y. Sire and E. Valdinoci, Fractional laplacian phase transitions and
boundary reactions: a geometric inequality and a symmetry result,
{\it J. Funct. Anal. 256}, 1842-1864 (2009).

\bibitem {St}G. Stampacchia, Some limit cases of $L^p$-estimates for solutions of second order elliptic equations, {\it Comm.
Pure Appl. Math. 16}, 505-510 (1963).

\bibitem {Ste} E. Stein, Singular Integrals and Differentiability Properties of Functions,
{\it Princeton University Press} (1970).

\bibitem {T} L. Tartar, Sur un lemme d'\'{e}quivalence utilis\'{e} en analyse num\'{e}rique, {\it Calcolo  24}, 129-140 (1987).

\bibitem {LV1}  L. V\'{e}ron, Singular solutions of some nonlinear elliptic equations,
{\it Nonlinear Anal. T., M. \& A.  5},  225-242, 1981.


\bibitem {V}  L. V\'{e}ron, Elliptic equations involving Measures,
 Stationary Partial Differential equations,
{\it Vol. I, 593-712, Handb. Differ. Equ., North-Holland, Amsterdam}
(2004).

\bibitem {V1}  L. V\'{e}ron, Existence and Stability of Solutions of General Semilinear
Elliptic Equations with Measure Data, {\it Advanced Nonlinear
Studies 13}, 447-460  (2013).

\end{thebibliography}
\end{document}